\documentclass[12pt]{amsart}
\usepackage{amssymb,amsmath}
\makeatletter
\def\@cite#1#2{{\m@th\upshape\bfseries%
[{#1\if@tempswa{\m@th\upshape\mdseries, #2}\fi}]}}
\makeatother
\theoremstyle{plain}
\newtheorem{thm}{Theorem}[section]
\newtheorem{lem}[thm]{Lemma}
\newtheorem{cor}[thm]{Corollary}
\newtheorem{prop}[thm]{Proposition}
\theoremstyle{definition}
\newtheorem{rem}[thm]{Remark}

\newtheorem{defn}[thm]{Definition}
\newtheorem{note}[thm]{Note}

\newtheorem{assumption}[thm]{Assumption}
\newtheorem{eg}[thm]{Example}
\newtheorem{egs}[thm]{Examples}

\newcommand{\Prf}{\noindent\textbf{Proof.\ }}
\newcommand{\bx}{\strut\hfill$\blacksquare$\medbreak}

\newcommand{\ca}{\mathrm{C}^*}

\newcommand{\ol}{\overline}

\newcommand{\td}{\widetilde}

\DeclareMathOperator*{\wotalg}{\textsc{wot}--{\rm Alg}}

\DeclareMathOperator*{\sotlim}{\textsc{sot}--lim}

\newcommand{\wot}{\textsc{wot}}

\newcommand{\bbA}{{\mathbb{A}}}

\newcommand{\bbC}{{\mathbb{C}}}

\newcommand{\bbF}{{\mathbb{F}}}

\newcommand{\bbQ}{{\mathbb{Q}}}

 \newcommand{\A}{{\mathcal{A}}}
 \newcommand{\B}{{\mathcal{B}}}
 \newcommand{\C}{{\mathcal{C}}}

\renewcommand{\H}{{\mathcal{H}}}
 
 \newcommand{\J}{{\mathcal{J}}}
 \newcommand{\K}{{\mathcal{K}}}  

 \newcommand{\M}{{\mathcal{M}}}
 \newcommand{\N}{{\mathcal{N}}}

\renewcommand{\S}{{\mathcal{S}}}

 \newcommand{\W}{{\mathcal{W}}}

\newcommand{\upchi}{{\raise.35ex\hbox{$\chi$}}}
\newcommand{\fA}{{\mathfrak{A}}}

\newcommand{\fC}{{\mathfrak{C}}}

\newcommand{\fL}{{\mathfrak{L}}}

\newcommand{\fR}{{\mathfrak{R}}}
\newcommand{\fS}{{\mathfrak{S}}}

\newcommand{\qand}{\quad\text{and}\quad}

\newcommand{\qor}{\quad\text{or}\quad}
\newcommand{\qfor}{\quad\text{for}\quad}

\newcommand{\qwith}{\quad\text{with}\quad}
\newcommand{\qonlyif}{\quad\text{only if}\quad}

\newcommand{\qiff}{\quad\text{if and only if}\quad}
\newcommand{\Alg}{\operatorname{Alg}}

\newcommand{\dist}{\operatorname{dist}}

\newcommand{\Id}{{\operatorname{Id}}}

\newcommand{\Lat}{\operatorname{Lat}}

\newcommand{\ran}{\operatorname{Ran}}
\newcommand{\rank}{\operatorname{rank}}

\newcommand{\spn}{\operatorname{span}}

\newcommand{\tr}{\operatorname{tr}}

\newcommand{\fngee}{\bbF^+\!(G)}

\newcommand{\sumin}{\sum_{i=1}^n}

\newcommand{\bofh}{\B(\H)}

\newcommand{\rows}{(S_1, \ldots, S_n)}

\newcommand{\fnplus}{{\bbF}_n^+}

\newcommand{\flgee}{\fL_G}
\newcommand{\frgee}{\fR_G}
\newcommand{\fgeeplus}{\bbF^+(G)}

\begin{document}

\title[Ideal Structure in Free Semigroupoid Algebras]%
{Ideal Structure In Free Semigroupoid Algebras From Directed 
Graphs}
%
\author[M.T.Jury and D.W.Kribs]{Michael~T.~Jury and 
David~W.~Kribs} 
\thanks{{\it key words and phrases.} Hilbert space, Fock space, 
directed graph, partial isometry, nonselfadjoint operator algebra, 
partly free algebra, Wold 
Decomposition, distance formula, Carath\'{e}odory Theorem.} 
\address{Department of Mathematics, Purdue University, West Lafayette, 
IN, USA 47907}
\email{jury@math.purdue.edu}
\address{Department of Mathematics and Statistics, University of Guelph, Guelph, ON,
CANADA N1G 2W1}
\email{dkribs@uoguelph.ca}

\date{}
\begin{abstract}
A {\it free semigroupoid algebra} is the weak operator topology closed 
algebra generated by the left regular
representation of a directed graph. 
We establish lattice 
isomorphisms between ideals and invariant subspaces, and this leads to a 
complete description of the $\wot$-closed ideal structure for these algebras. We 
prove a distance formula to ideals, and this gives an appropriate version of 
the Carath\'{e}odory interpolation theorem. Our 
analysis rests on an investigation of predual properties, specifically 
the $\bbA_n$ properties for linear functionals, together with a general 
Wold Decomposition for $n$-tuples of partial isometries. 
A number of our proofs unify proofs for subclasses appearing in the 
literature. 
\end{abstract} 
\maketitle



In \cite{KP2} and \cite{KP1}, the second author and Stephen Power began studying a
class of operator algebras called {\it free semigroupoid algebras}. These are the $\wot$-closed
(nonselfadjoint) algebras $\flgee$ generated by the left regular representations of directed graphs
$G$. Earlier work of Muhly and Solel \cite{M1,MS1} considered the norm closed algebras
generated by these representations in the finite graph case; they called them {\it quiver
algebras}. 
In the case of single vertex graphs, the $\flgee$ obtained include the classical analytic Toeplitz
algebra $H^\infty$
\cite{Douglas_text,Hoffman_text,Sar1} and the noncommutative analytic Toeplitz algebras
$\fL_n$
studied by Arias, Popescu, Davidson, Pitts,  and others 
\cite{AP,DKP,DP1,DP2,Kfactor,Pop_fact,Pop_beur}.

In this paper, we consider algebraic structure-type problems for the 
algebras $\flgee$. 
In particular, we derive a complete description of the $\wot$-closed ideal
structure; for instance, there is a lattice isomorphism between right ideals and
invariant subspaces of the commutant $\flgee^\prime = \frgee$ \cite{DP2}. Furthermore, we
prove a distance formula to ideals in these algebras \cite{AP2,DP3,McCull}.
This yields a version of the Carath\'{e}odory interpolation theorem \cite{AP2,DP3} for $\flgee$. 
A valuable tool in our analysis is 
a general Wold Decomposition \cite{MS1,Pop_diln}, which we establish for $n$-tuples of
partial
isometries with initial and final projections satisfying natural conditions. This
leads to information on predual properties for $\flgee$. We prove $\flgee$ satisfies property
$\bbA_1$ \cite{BFP,DP2}; that is, every weak-$*$ continuous linear functional may
be realized as a vector functional. A number of our proofs for general 
$\flgee$ unify the proofs for the  
special cases of  $\fL_n$ and $H^\infty$, which were previously 
established by different means. 

The first section contains a brief review of the notation associated with these algebras. Our
attention will be focused on the cases when the directed graph $G$ has `no 
sinks'; that is, every vertex in $G$ is the initial vertex for some directed edge. We 
include a list of some examples generated by simple graphs. We begin the analysis proper in the
second
section, with a Wold Decomposition for $n$-tuples of partial isometries.   
This leads
into the topic of the third section; an investigation into the basic 
properties of linear functionals on
$\flgee$. In particular, we show the ampliation algebras 
$\flgee^{(n)}$  have the factorization property $\bbA_n$. 
In the subsequent section, we prove  the subclass of algebras with {\it partly free} 
commutant discovered in \cite{KP2,KP1} are precisely those $\flgee$ which 
satisfy the stronger factorization property   $\bbA_{\aleph_0}$, when an initial restriction is
made on the graph.  
Using property $\bbA_1$ for $\flgee$ and the Beurling Theorem from \cite{KP1}, we establish
complete lattice isomorphisms between ideals and invariant subspaces of $\flgee$ in the fifth 
section. This allows us to describe, for example, the $\wot$-closure of the commutator ideal, 
and
precisely when $\wot$-closed ideals are finitely generated. In the penultimate section we prove a
completely isometric 
distance formula to ideals of $\flgee$, and we apply this in special cases to obtain a
Carath\'{e}odory Theorem in the final section.

\section{Free Semigroupoid Algebras}\label{S:freesemigroupoid}

Let $G$ be a countable (finite or countably infinite) directed graph with edge set $E(G)$ and 
vertex set $V(G) $. Let $\fngee$ be the free semigroupoid determined by 
$G$; that is, $\fngee$ consists of the vertices which act as units, written as $\{k \}_{k\geq 1}$, 
and   allowable finite paths in
$G$, with the
natural operations of 
concatenation of allowable paths. Given a path $w= e_{i_m}\cdots e_{i_1}$ in $\fngee$, an
allowable product of edges $e_{i_j}$ in $E(G)$, we
write $w=k_2 wk_1$ when the initial and final vertices of $w$ are, respectively, $k_1$ 
and $k_2$. Further, by $|w|$ we mean the number of directed edges which determine the path
$w$. 

\begin{assumption}
For our purposes, it is natural to restrict attention to directed graphs 
$G$ with `no sinks'; that is, every vertex is the initial vertex for some 
directed edge. 
\end{assumption}

Let $\H_G  = \ell^2(\fngee)$ be the Hilbert space with 
orthonormal basis $\{ \xi_w : w\in \fngee\}$ indexed by elements of $\fngee$. For each edge
$e\in E(G)$ and vertex $k \in V(G)$, define partial isometries and projections on $\H_G$ by: 
\[
L_e\xi_w = \left\{ \begin{array}{cl}
\xi_{ew} & \mbox{if $ew\in\fngee$} \\
0 & \mbox{otherwise}
\end{array}\right.
\]
and
\[
L_k \xi_w = \left\{ \begin{array}{cl}
\xi_{kw}=\xi_w & \mbox{if $w = kw\in\fngee$} \\
0 & \mbox{otherwise}
\end{array}\right.
\]

These operators may be regarded as `partial creation operators' 
acting on a generalized
Fock space Hilbert space. There is an equivalent tree perspective, discussed in
\cite{KP1}, which gives an appealing visual interpretation of the actions of 
these operators. The vectors $\{ \xi_k : k\in V(G)\}$ are called the {\it vacuum vectors}.

The family $\{L_e,L_k\}$ also arises through the left regular representation
$\lambda: \fngee \rightarrow \B(\H_G)$, with $\lambda(e) = L_e$ and $\lambda(k) =
L_k$. 
The associated {\it free semigroupoid algebra} is the weak operator 
topology closed algebra generated by this family; 
\begin{eqnarray*}
\flgee &=& \wotalg\,\, \{ L_e,L_k :e\in E(G), k\in V(G) \} \\
&=& \wotalg \,\,\{ \lambda(w) : w\in\fngee \}.
\end{eqnarray*}
These algebras were the focus of analysis by the second author and Power in \cite{KP2,KP1}. 
In the case of finite graphs, Muhly and Solel \cite{M1,MS1} considered the 
norm closed algebras $\A_G$ generated by such a family, calling them {\it 
quiver algebras}. 

There is an analogous right regular representation $\rho:\fngee\rightarrow
\B(\H_G)$, which yields partial isometries $\rho(w) \equiv R_{w^\prime}$ for
$w\in\fngee$ acting on $\H_G$ by the equations $R_{w^\prime}\xi_v = \xi_{vw}$, where
$w^\prime$ is the word $w$ in reverse order, with similar conventions. 
Observe that
$R_{v^\prime} L_w = L_w R_{v^\prime}$ for all $v,w\in\bbF^+(G)$. In fact, the algebra  
\begin{eqnarray*}
\frgee &=& \wotalg \,\, \{R_e,R_k :e\in E(G), k\in V(G) \} \\
&=& \wotalg \,\, \{ \rho(w) : w\in\fngee \}
\end{eqnarray*}
coincides with the commutant $\flgee^\prime = \frgee$. The commutant is also unitarily
equivalent to the algebra $\fL_{G^t}$, where $G^t$ is obtained from $G$ simply by reversing
directions of all edges.  Elements of $\flgee$ have Fourier
expansions:  If $A\in\flgee$ and $k \in V(G)$, then
$A\xi_k = \sum_{w=wk} a_w \xi_w$ for some scalars $a_w\in\bbC$, and the 
Cesaro-type sums 
\[
\Sigma_k(A) = \sum_{ |w|<k} \Big( 1 -
\frac{|w|}{k}\Big) a_w L_w.
\]
converge in the strong operator topology to $A$. These results are contained in \cite{KP1}. We
write $A\sim
\sum_{w\in\fngee} a_w L_w$ as a notational convenience.  We shall also put  $P_k =
L_k$ and $Q_k = R_k$ for the projections determined by vertices $k \in V(G)$.

A valuable tool in our analysis is 
the Beurling-type  invariant subspace theorem for $\flgee$ proved in 
\cite{KP1}. 
A non-zero subspace $\W$ of $\H_G$ is {\it wandering} for $\flgee$
if the subspaces $L_w\W$ are pairwise orthogonal for distinct $w$ in
$\fngee$.
Observe that every wandering subspace $\W$ generates an
$\flgee$-invariant subspace given by 
\[
\flgee [\W] = \sum_{w\in\fngee} \oplus L_w \W.
\]
Every $\flgee$-wandering vector $\zeta$ generates the cyclic invariant subspace
$\flgee[\zeta]$. The subspace  $\flgee[\zeta]$ is  {\it minimal cyclic} if $L_k \zeta = P_k \zeta =
\zeta$ for some vertex $k\in V(G)$. Notice that if $\zeta$ is a wandering 
vector, then each
vector $P_k \zeta$ which is non-zero is wandering as well. 
The following result was proved in \cite{KP1}. It is a generalization of Beurling's classical
theorem \cite{Beur} for $H^\infty$ and a corresponding result for free semigroup algebras
$\fL_n$ \cite{AP,DP1}. 

\begin{thm}\label{beurlingthm}
Every invariant subspace of $\flgee$ is generated by a wandering subspace,
and is the direct sum of minimal cyclic subspaces generated by wandering
vectors. Every minimal cyclic invariant subspace generated by a
wandering vector is the range of a partial isometry in $\frgee$, and the choice 
of partial isometry is unique up to a scalar multiple. 
\end{thm}

In fact, more was proved in \cite{KP1}. The partial isometries in $\frgee$  have a standard form,
and their initial projections are sums of projections amongst $\{Q_k: k\in V(G)\}$. 
Given a minimal cyclic subspace $\flgee[\zeta]$ with $P_k \zeta = \zeta$, a partial isometry
$R_\zeta$ in $\frgee$
which satisfies $\flgee[\zeta] = R_\zeta \H_G$ is defined by 
$R_\zeta\xi_w= L_w \zeta$ for $w$ in $\fngee$, and the initial projection satisfies $R_\zeta^*
R_\zeta = Q_k$. Further, any partial isometry in $\frgee$ with range space  $\flgee[\zeta]$ is a
scalar multiple of $R_\zeta$.

We finish this section by setting aside a number of examples generated by simple graphs.

\begin{egs}\label{examples}
$(i)$ The algebra generated by the graph with a single vertex and single loop edge is unitarily
equivalent to the classical analytic Toeplitz algebra $H^\infty$
\cite{Douglas_text,Hoffman_text,Sar1}. 

$(ii)$ The noncommutative analytic Toeplitz algebras $\fL_n$, $n\geq 2$
\cite{AP,AP2,DKP,DP1,DP2,DP3,Kfactor,Pop_fact,Pop_beur}, arise from the graphs with a
single
vertex and $n$ distinct loop edges. 

$(iii)$ The cycle algebras $\fL_{C_n}$ discussed in \cite{KP1} are generated by the graph
$C_n$ with $2 \leq n< \infty$ vertices, and edges corresponding to directions $\big\{ (2,1),
\ldots, (n,n-1), (1,n)\big\}$. These algebras may be represented as matrix function algebras. 

$(iv)$ The algebra $\fL_{C_\infty}$ \cite{KP2,KP1} generated by the infinite graph analogue of
the cycle graphs is determined by the graph $C_\infty$ with vertices $\{ k \}_{k\geq 1}$ and
directed edges $(k+1,k)$. 

$(v)$  A non-discrete example discussed in \cite{KP2} is given by the 
graph $Q$  consisting of
vertices $\{ q\}_{q \in\bbQ}$ indexed by the rational numbers, and directed edges $e_{qp} = q
e_{qp} p$
whenever $p \leq q$. Notice that $Q^t$ is 
graph isomorphic to $Q$, and thus the algebra $\fL_Q \simeq \fL_{Q^t} \simeq \fL_Q^\prime$
is unitarily equivalent to its commutant. 
\end{egs}


\section{Wold Decomposition}\label{S:wold}

In this section, we establish a Wold
Decomposition for $n$-tuples of (non-zero) operators $S = \rows$ which act on a common
Hilbert
space $\H$  and  satisfy the following relations denoted by $(\dagger)$: 
\begin{itemize}
\item[$(1)$] 
$
\sumin S_iS_i^* \leq I.
$  
\end{itemize}
\begin{itemize}
\item[$(2)$] For $1 \leq i \leq n$, 
\[
(S_i^* S_i)^2 = S_i^* S_i .
\] 
\item[$(3)$] For $1\leq i,j \leq n$, 
\[
(S_i^*S_i)(S_j^* S_j) = 0 \qor S_i^* S_i = S_j^* S_j .
\]
\item[$(4)$] 
For $1\leq i \leq n$, there is a $j$ such that 
\[
S_i S_i^* \leq S_j^* S_j.   
\]
\item[$(5)$]
The distinct elements $\{P_k \}_{k\in\S}$ from the set $\{ 
S_i^* S_i  : 1\leq i \leq n \}$ satisfy
\[ 
\sum_{k\in\S} P_k = I .
\]
\end{itemize}

\begin{note}
We shall behave as though $n$ is finite in this section; there are obvious modifications which
can be made in the infinite case. 
The first four conditions say the $S_i$ are partial isometries with 
pairwise orthogonal ranges, with initial projections either orthogonal or 
equal, and with range projections supported on a (unique) initial 
projection. The last condition, which is redundant when equality is 
achieved in the first condition, means no non-zero vector is 
annihilated by all the $S_i$. In terms of the directed graph connection we 
are about to make, this means the associated directed graphs have no 
sinks. 
\end{note}

\begin{defn}
Given a positive integer $n
\geq 2$, we write $\fnplus$ for the (non-unital) free semigroup on $n$ noncommuting letters.
If $G$ is a directed graph with $n$ edges, we let $\fnplus (G)$ denote the set of all finite 
words  in the edges of $G$. Observe that $\fnplus (G)$ contains the set of finite paths
$\fgeeplus \setminus V(G)$ as a subset. 
\end{defn}

If
$w =i_m\cdots i_1$ belongs to $\fnplus$, it is convenient in this section to  let $w(S)$ be the
operator product $w(S) =
S_{i_m} \cdots S_{i_1}$. When $S = \rows$ satisfies $(\dagger)$, this is a partial isometry with
initial projection $w(S)^* w(S)$ equal to $S_{i_1}^* S_{i_1}$ or $0$, and final projection
$w(S) w(S)^* $ supported on $S_{i_m} S_{i_m}^*$. 
Let $\fS_S$ be  the weak operator topology closed algebra generated by such an
$n$-tuple and its initial projections;
\[
\fS_S = \wot - \Alg \big\{ S_1, \ldots, S_n, S_1^*S_1, \ldots, S_n^* S_n \big\}. 
\]

\begin{defn}\label{wandering}
A subspace $\W$ of $\H$ is {\it wandering} for $S = \rows$ satisfying $(\dagger)$ if the
subspaces $w(S) \W$ are
pairwise orthogonal for distinct words $w$ in $\fnplus$. Observe that a given partial isometry
$w(S)$ may be equal to zero here. Every wandering subspace generates an invariant subspace by 
\[
\fS_S[\W] = \sum_{w\in\fnplus} \oplus w(S) \W . 
\]
\end{defn}

The following is the pure part of the Wold Decomposition. We set it aside since it is
the precise form we require in this paper. 

\begin{lem}\label{purewold}
Let $S = \rows$ be operators on $\H$ satisfying $(\dagger)$. The subspace $\W = \ran ( I -
\sumin S_i S_i^*)$ is wandering for $S$. Suppose that  
\[
\H = \sum_{w\in\fnplus} \oplus w(S) \W. 
\]
Let $\{ P_k \}_{k\in\S}$ be the distinct projections from the set $\{ S_i^* S_i : 1\leq i \leq n\}$.
Then 
\[
P_k \W \subseteq \W = \sum_{k\in\S} \oplus P_k \W . 
\]

Let $G$ be the directed graph (with no sinks)  with vertex set $V(G) 
\equiv \S$, and $n$ directed
edges $\{ e \in E(G) \}$ where the number of edges from vertex $k$ to vertex $l$ is given by 
the cardinality of the set 
\[
\{ S_i \,\,:\,\, S_i^* S_i =P_k \qand S_i S_i^* \leq P_l \}. 
\]
Given $k\in \S$, let $\alpha_k = \dim P_k \W$ and let $\{ S_i^{(k)}\} = \{ S_i : S_i^* S_i = P_k
\}$. The sets $\{ S_i^{(k)}\}$ and $\{ L_e : e=ek \}$ have the same
cardinality for each $k\in \S = V(G)$.  Let $\S_0 = \{ k\in V(G) : \alpha_k\neq 0\}$. Then there
is a unitary 
\[
U : \H \longrightarrow \sum_{k\in \S_0} \oplus (Q_k \H_G)^{(\alpha_k)} 
\]
such that for $l \in \S$,
\[
\big\{ U S_i^{(l)} U^* \big\} = \left\{ \sum_{k\in \S_0}\oplus L_e^{(\alpha_k)}\Big|_{(Q_k
\H_G)^{(\alpha_k)}} : e=e l \right\},
\]
where $L_e^{(\alpha_k)}|_{(Q_k
\H_G)^{(\alpha_k)}}$ is the restriction of the ampliation $ L_e^{(\alpha_k)}$ to the
$\alpha_k$-fold direct sum of $Q_k \H_G$ with itself. 
\end{lem}

\Prf
As $S =( S_1, \ldots, S_n)$ are partial isometries with pairwise orthogonal ranges, the operator
$P = I
- \sumin S_i S_i^*$ is a projection and the subspace $\W = P \H$ is wandering. Notice that 
\[
( S_i^* S_i) ( S_j S_j^*) = ( S_j S_j^*)  ( S_i^* S_i) 
= \left\{ \begin{array}{cl}
S_j S_j^* & \mbox{if $S_j S_j^* \leq S_i^* S_i$} \\ 
0 & \mbox{otherwise}
\end{array}\right. 
\]
Thus the initial projections $S_i^* S_i$ commute with $P$, and hence $P_k \W$ is contained in
$\W$ for $k\in\S$. Further, $\W = \sum_{k\in\S} \oplus P_k \W$ by 
condition $(5)$ in $(\dagger)$; the no-sink
condition. 
The sets $\{ S_i^{(k)} \}$ and $\{ L_e : e=ek,\,\, e\in E(G) \}$ have the same cardinality from
the
definition of $G$. It is convenient to re-label $\{ S_i^{(k)} \}$ as  $\{ S_e : e=ek, e\in E(G) \}$
for the rest of the proof.

It remains to construct a unitary $U$ which intertwines $\{S_e : e\in E(G) \}$ with the
appropriate sums of restricted ampliations of $\{L_e : e\in E(G) \}$. 
First note that there is a natural bijective correspondence between $\fnplus (G)$ and $\fnplus$
induced by the re-labelling of $\{S_1, \ldots, S_n \}$ as $\{S_e : e\in E(G) \}$. The non-trivial
paths in the free
semigroupoid, $\fgeeplus\setminus V(G)$, may be regarded as a subset of 
$\fnplus (G)$. The key
point is that,
under this identification, the set of all words $w$ in $\fnplus (G)$ for which the partial isometry
$w(S)$ is non-zero is precisely the set $\fgeeplus\setminus V(G)$. Indeed,
given a formal product $w = e_k \cdots e_1$ in $\fnplus (G)$, the definition of $G$ and
properties $(\dagger)$ for $S = \rows$ show that 
\begin{eqnarray*}
w(S) \neq 0 & \qiff & S_{e_j} S_{e_j}^* \leq  S_{e_{j+1}}^* S_{e_{j+1}} \qfor 1\leq j < k \\ 
& \qiff & w = e_k \cdots e_1 \in \fgeeplus
\end{eqnarray*}

Given $k$ in $\S_0$, so that $\alpha_k = \dim P_k\W \neq 0$, for $1\leq j \leq \alpha_k$ let 
\[
\Big\{ \xi_w^{(j)} \equiv w(L) \xi_k^{(j)} : w\in \fgeeplus, \, w=wk\Big\}
\]
be the standard orthonormal
basis for the $j$th copy of the subspace $Q_k \H_G$ in the $\alpha_k$-fold direct sum
$(Q_k\H_G)^{(\alpha_k)}$. Next, for each $k$ in $\S_0$ choose an orthonormal basis $\{
\eta_k^{(j)} : 1 \leq j \leq \alpha_k \}$  for the (non-zero) subspace $P_k \W$. By hypothesis we
have 
\[
\H = \spn \Big\{ w(S) \eta_k^{(j)} : w\in \fnplus(G), \, k\in\S_0, \, 1\leq j \leq \alpha_k \Big\}.
\]
Moreover, as $\W = \sum_{k\in\S_0}\oplus P_k \W$ is wandering for $\S$, the non-zero vectors
in
this spanning set form an orthonormal basis for $\H$. But we observed above that $w$ in
$\fnplus(G) \setminus \fgeeplus$ implies $w(S) = 0$. Further, by the 
properties $(\dagger)$ we also have $w(S) \eta_k^{(j)} \neq 0$ precisely 
when $w=wk$ inside $\fgeeplus$. Hence the spanning set for $\H$ may be
restricted to require $w$ in $\fgeeplus\setminus V(G)$. Thus, it follows that we may define a
unitary operator 
\begin{eqnarray*}
U : \H   \longrightarrow \sum_{k\in\S_0} \oplus (Q_k \H_G)^{(\alpha_k)}
\end{eqnarray*}
by intertwining these bases; 
\[
\left\{ \begin{array}{ccc} 
U\big( \eta_k^{(j)} \big) &=& \xi_k^{(j)} \\
 & & \\  
U\big( w(S) \eta_k^{(j)}\big) & = & w(L) \xi_k^{(j)} \,\, (= \xi_w^{(j)}) 
\end{array}\right. 
\]
where $k\in\S_0$, $1\leq j \leq \alpha_k$, and $w=wk\in \fgeeplus\setminus V(G)$. 
Evidently, this unitary jointly intertwines the operators $\{S_i \}$ with the ampliations 
$\big\{ \sum_{k\in\S_0} \oplus L_e^{(\alpha_k)}|_{(Q_k \H_G)^{(\alpha_k)}}\big\}$ in the
desired manner, and this completes the proof.
\bx

It follows that the algebra $\fS_S$ is unitarily equivalent to a corresponding `weighted
ampliation' of $\flgee$, with weights $\alpha_k$. We may regard the dimension $\alpha_k =
\dim P_k \W$ as the {\it pure multiplicity over vertex} $k$ in the Wold Decomposition for $S$.
Observe that we are allowing for $\alpha_k = 0$; this means that no copy of the `tree
component' subspace $Q_k \H_G$ appears in the range of the unitary $U$. 
This situation arises, for instance, if the operators $L = (L_e)_{e\in E(G)}$ themselves are
restricted to the direct sum of a  subset of these reducing 
subspaces $\{Q_k \H_G \}_{k\in V(G)}$. More generally, there is an interesting consequence of
the lemma in the case that $S$ is a joint restriction of some 
$L = (L_e)_{e\in E(G)}$ to a cyclic
invariant subspace. 
A version of the following result for the ampliations of $\flgee$ will 
play a key role in Sections \ref{S:predual} and \ref{S:distance}. 

\begin{cor}\label{cyclicform}
Let $\M = \ol{\flgee \xi}$ be a cyclic invariant subspace for $\flgee$. Let $\S = \{ k\in V(G) :
P_k \xi \neq 0 \}$. Then the operators $L_e|_{\M}$ are jointly unitarily equivalent to the
operators $L_e|_{\sum_{k\in\S} \oplus (Q_k \H_G)}$. In particular, the corresponding
restrictions of $\flgee$ are unitarily equivalent; 
\begin{eqnarray}\label{cyclic}
\flgee|_\M \simeq \flgee|_{\sum_{k\in\S} \oplus (Q_k \H_G)}. 
\end{eqnarray}
\end{cor}

\Prf
Assume first that all the restrictions $S_e = L_e|_{\M}$ are non-zero. Then we may apply
Lemma~\ref{purewold} 
to the tuple
$(S_e)_{e\in E(G)}$. The wandering subspace here is 
\[
\W = \ran \Big( I_\M - \sum_{e\in E(G)} S_e S_e^* \Big) = \dim \left( \M 
\ominus \Big( \sum_{e\in
E(G)} \oplus L_e \M \Big)\right). 
\]
But in this case $\W$ is spanned by $\{ P_k P \xi = P P_k \xi : k \in V(G) 
\}$ where $P = I - \sum_e
L_e L_e^*$. It 
follows that the subspaces $P_k \W \subseteq \W$ are at most one 
dimensional, and that
$P_k\W \neq 0$ if and only if $P_k \xi \neq 0$. Hence the result follows 
from Lemma~\ref{purewold}.

In the general case, Lemma~\ref{purewold} would produce a joint unitary equivalence between
the non-zero operators amongst $\{ S_e : e\in E(G)\}$ and the creation operators from a
subgraph of $G$. (The subgraph obtained would consist of the saturations at all vertices in $\S =
\{ k\in V(G) : P_k\xi\neq 0 \}$.) But it is easy to see that $S_e = 0$ if and only if $e$ is not part
of a path starting at some vertex in $\S$, and this also corresponds to the case that
$L_e|_{\sum_{k\in\S} \oplus (Q_k\H_G)^{(\alpha_k)}} = 0$. Thus, the unitary equivalence
(\ref{cyclic}) holds in all situations. 
\bx

We now turn our attention to the general case. The families of partial 
isometries $L = (L_e)_{e\in
E(G)}$ provide the models for pure partial isometries in the Wold 
Decomposition through their weighted
ampliations as in Lemma~\ref{purewold}. On the other hand, the 
`coisometric' component of the decomposition is characterized by 
determining a representation of a  Cuntz-Krieger directed graph 
$\ca$-algebra \cite{BHRS,Kumjian1,Kumjian2}. 

\begin{defn}
Let $S = \rows$ satisfy $(\dagger)$ on $\H$. Then 
\begin{itemize}
\item[$(i)$]
$S$ is {\it fully coisometric} if equality is achieved in condition $(1)$ 
of $(\dagger)$; 
\[
SS^* = 
\left[
\begin{matrix}
S_1 & \cdots & S_n 
\end{matrix} \right] \,\,
\left[ \begin{matrix}
S_1^* \\ 
\vdots \\ 
S_n^* 
\end{matrix} \right]
= \sumin S_i S_i^* = I. 
\]
\item[$(ii)$] 
$S$ is {\it pure} if $\H$ is equal to the $S$-invariant subspace generated by the wandering
subspace $\W = \ran \Big( I - \sumin S_iS_i^* \Big)$; 
\[
\H = \sum_{w\in\fnplus} \oplus w(S) \W. 
\]
\end{itemize}
\end{defn}

Every countable directed graph $G$ with $n$ edges can be seen to determine a fully coisometric
$n$-tuple
$S= (S_e)_{e\in E(G)}$, where the particular $(\dagger)$ relations between 
the $S_e$ are governed
by the directed graph as in the statement of Lemma~\ref{purewold}. Conversely,  every
fully coisometric $n$-tuple here can be shown to be jointly unitarily equivalent to a fully
coisometric $n$-tuple of the form $S = (S_e)_{e\in E(G)}$, where the 
directed graph $G$ is again
explicitly obtained as in Lemma~\ref{purewold}. (See work of Brenken  
\cite{Brenken} for discussions on this topic.) Furthermore, 
we have shown 
in Lemma~\ref{purewold} how pure $n$-tuples are completely determined by 
tuples $L = (L_e)_{e\in
E(G)}$. Thus we may prove the following. 

\begin{thm}\label{wold}
Let $S=\rows$ be operators on $\H$ satisfying $(\dagger)$. Then
$S_1, \ldots, S_n$ are jointly unitarily equivalent to the direct sum of a pure $n$-tuple and a
fully coisometric $n$-tuple which both satisfy $(\dagger)$. In other words, there is a directed
graph $G$ with $n$ edges and a unitary $U$ such that  $US_1U^*,\ldots, US_nU^* $ are of the
form 
\begin{eqnarray}\label{woldid}
\,\,\,\,\,\,\,\,\, \big\{ US_iU^* \big\}_{i=1}^n = 
\left\{ S_e \oplus \Big( \sum_{k\in \S_0}\oplus L_e^{(\alpha_k)}\Big|_{(Q_k
\H_G)^{(\alpha_k)}} \Big) : e\in E(G) \right\}. 
\end{eqnarray}

Let $\H_p = \sum_{w\in\fnplus} \oplus w(S) \W$ where $\W = \ran \big( I - \sumin S_i S_i^*
\big)$ and let $\H_c = ( \H_p)^\perp$. 
The subspaces $\H_c$ and $\H_p$ reduce $S = \rows$, and the restrictions $S_i|_{\H_c}$ and 
$S_i|_{\H_p}$ determine the joint unitary equivalence in $(\ref{woldid})$. 

This decomposition is unique in the sense that if $\K$ is a
subspace of $\H$ which reduces $S = \rows$, and if the restrictions $\{ S_i|_\K : 1\leq i \leq n\}$
are pure, respectively fully coisometric, then $\K \subseteq \H_p$, respectively $\K \subseteq
\H_c$. 
\end{thm}

\Prf
Since each $S_i^* \W = \{0\}$, it is clear that $\H_p$ reduces $S = \rows$. Thus, the restrictions
$S_i|_{\H_p}$ form a pure $n$-tuple and Lemma~\ref{purewold} shows they are jointly
unitarily equivalent to an $n$-tuple of the desired form. On the other hand, the operators $V_i =
S_i|_{\H_c}$ satisfy 
\[
I_{\H_c} = \sumin V_i V_i^* = \sumin S_iS_i^*|_{\H_c}.
\]
Indeed, if $P
= I - \sumin S_iS_i^*$, then $P\xi$ belongs to $\W \subseteq \H_p$ for all $\xi\in\H$. But when
$\xi\in\H_c$ we have $P\xi\in\H_c$ since $\H_c$ reduces $S$, and thus $P\xi \in \H_c \cap \W
\subseteq \H_c \cap \H_p = \{ 0 \}$, so that $P\xi = 0$ as claimed.  
Hence it follows that $V = (V_1, \ldots, V_n)$ satisfies $(\dagger)$ and is fully coisometric, and
by our remarks preceding the theorem, this $n$-tuple is determined by the same directed graph
$G$ as the pure part of $S = \rows$. 

To verify the uniqueness statement we shall prove more. Let $P_c$ and
$P_p$ be, respectively, the projections of $\H$ onto $\H_c$ and $\H_p$.  
Let $Q\in\bofh$ be a projection such that $Q\H$ reduces $S = \rows$. We
claim that $QP_c = P_c Q $ and $QP_p = P_p Q$; in other words, $Q\H = P_c
(Q \H) \oplus P_p (Q\H)$ contains both subspaces $P_c (Q\H)$ and $P_p 
(Q\H)$.  
To see this, first let $\{ X_k : k\geq 1\}$ be the projections 
\[ X_k =
\sum_{w\in\fnplus ; \, |w|=k} w(S) w(S)^* = \Phi^k (I), 
\] 
where $\Phi(A)
= \sumin S_i A S_i^*$. As $\H_c$ reduces $S$ and $\Phi(P_c)  = P_c$, we
have the restrictions $X_k|_{\H_c} = I_{\H_c}$ for $k\geq 1$.  
Furthermore, since $\H_p = \sum_{w\in\fnplus} \oplus w(S) \W$ and $w(S)^*
\W = \{ 0\}$ for $w\in \fnplus$, it is evident that
the strong operator topology limit $\sotlim_{k\rightarrow\infty} 
X_k|_{\H_p} = 0$.

Now let $Q\xi = \xi = P_c\xi + P_p\xi$ belong to $Q\H$. It suffices to 
show that both $P_c\xi$ and $P_p\xi$ belong to $Q\H$. But since $Q\H$ is 
reducing for $S$, we have 
\[
P_c \xi = \lim_{k\rightarrow\infty} X_k P_c\xi = \lim_{k\rightarrow\infty} 
X_k (P_c\xi + P_p\xi) = \lim_{k\rightarrow\infty} X_k\xi  \in Q\H.
\]
This establishes the claim and finishes the proof. 
\bx

\begin{rem}              
This decomposition theorem plays a role in forthcoming work of the authors on dilation theory
\cite{JK2}. 
There are a number of modern generalizations of Wold's classical theorem which appear in the
literature (for instance see \cite{JSW,MS1,Pop_diln}). A 
special case of Theorem~\ref{wold} is discussed in 
\cite{KP1}. In its general form, Theorem~\ref{wold} is most closely
related to the Wold Decomposition of Muhly and Solel \cite{MS1} 
established for the more abstract setting of
representations of $\ca$-correspondences. In fact, from one point of view, 
Theorem~\ref{wold} can be thought of as an explicit identification, of 
the pure part in particular, of the components of their Wold Decomposition 
in some very concrete cases. 
\end{rem}

\section{Predual Properties}\label{S:predual}

A $\wot$-closed algebra $\fA$ has property $\bbA_n$, $1\leq n 
\leq \aleph_0$, if for every $n\times
n$ matrix $[\varphi_{ij}]$ of weak-$\ast$-continuous linear functionals on $\fA$, 
there are vectors
$\{ \zeta_i, \eta_j :1\leq i,j \leq n\} $ with 
\begin{eqnarray}\label{predualan}
\varphi_{ij} (A) = (A\eta_j, \zeta_i) \qfor A\in\fA \qand 1\leq i,j \leq 
n. 
\end{eqnarray}
These notions are discussed in detail in \cite{BFP}.
Given $1 \leq n \leq \aleph_0$, recall that the $n$th {\it ampliation} of 
$\fA$ is the $\wot$-closed algebra $\fA^{(n)}$ generated by the $n$-fold 
direct sums $A^{(n)} = A \oplus \ldots \oplus A$ where $A\in \fA$. 
The infinite ampliation of $\fA$, which we write as $\fA^{(\infty)}$,  
satisfies property $\bbA_{\aleph_0}$ \cite{BFP}, and hence 
property  $\bbA_n$ for $n\geq 1$. 
We prove the following for $\flgee$ in the case that $G$ has no sinks. 

\begin{thm}\label{a1thm}
$\flgee^{(n)}$ has 
property $\bbA_n$ for $1 \leq n \leq \aleph_0$.
\end{thm}

\Prf
Fix $1 \leq n \leq \aleph_0$ and let $\fA = \flgee^{(n)}$ and $\H = \H_G^{(n)}$. 
Let $[\varphi_{ij}]$ be an $n\times n$ matrix of weak-$\ast$-continuous 
functionals on $\fA$. We may regard these as functionals 
$\td{\varphi}_{ij}$ acting on $\fA^{(\infty)}$ by $\td{\varphi}_{ij} 
(A^{(\infty)}) \equiv \varphi_{ij}(A)$. Since $\fA^{(\infty)}$ has 
property $\bbA_{\aleph_0}$, hence $\bbA_n$,  there are vectors $\{ x_i, 
y_j \}_{1\leq i,j \leq n}$ in $\H^{(\infty)}$ such that 
\[
\td{\varphi}_{ij} (A^{(\infty)}) = (A^{(\infty)} x_i ,y_j) = 
(A^{(\infty)} x_i ,P_\M y_j), 
\]
where $P_\M$ is the projection onto the $\fA^{(\infty)}$-invariant 
subspace 
\[
\M = \bigvee_{1 \leq i \leq n} \ol{\fA^{(\infty)} x_i}.
\] 
Clearly we may assume  $P_\M y_j = y_j$ for $1\leq j \leq n$. 

Let $\td{L_e}$ and $\td{P_k}$ be the restrictions of the  generators of
$\fA^{(\infty)}$ to this invariant subspace; that is,  
\[
\td{L_e} = 
(L_e^{(n)})^{(\infty)}|_\M \qand \td{P_k} = (P_k^{(n)})^{(\infty)}|_\M.
\] 
The wandering subspace here is  
\[
\W = \ran \Big( I_\M 
- \sum_{e\in E(G)} \td{L_e} \td{L_e}^*\Big)
= \M \ominus \sum_{e\in E(G)} \td{L_e} \M 
\]
Then, as in Corollary~\ref{cyclicform}, it follows that $\alpha_k \equiv \dim \td{P_k} \W \leq
n$ for $k$ in $V(G)$. Further, if we let $\S = \{ k\in  V(G) : \td{P_k} \W \neq 0 \}$,
Lemma~\ref{purewold} gives a unitary  
\[
U \, : \, \M \longrightarrow \sum_{k\in \S} \oplus (Q_k \H_G)^{(\alpha_k)} 
\hookrightarrow \H
\]
which may be defined so that 
\[
U \td{L_e} U^* = \sum_{k\in \S} \oplus L_e^{(\alpha_k)}|_{(Q_k 
\H_G)^{(\alpha_k)}}.
\]

Now, given $1 \leq i,j \leq n$, let $\td{x_i}$ be the vector in $\H$ 
defined by 
\[
P_{U\M} (\td{x_i}) = Ux_i \qand P_{\H\ominus U\M} (\td{x_i}) =0,
\]
and similarly define $\td{y_j}$ in terms of $U y_j$. Moreover, for each  
$A^{(\infty)}$ in $\fA^{(\infty)}$ let $A_\S = U A^{(\infty)} U^*$. 
Then given $A$ in $\fA = \flgee^{(n)}$ and $1\leq i,j \leq n$, we have 
\begin{eqnarray*}
\varphi_{ij} (A) = \td{\varphi}_{ij}(A^{(\infty)}) &=& (A^{(\infty)} x_i,y_j) \\ 
&=& (U^* A_\S U x_i, y_j) \\
&=& ( A_\S U x_i, U y_j) = (A \td{x_i}, \td{y_j}).  
\end{eqnarray*}
This shows that $\flgee^{(n)}$ has property $\bbA_n$. 
\bx

An immediate consequence of the fact that $\flgee$ satisfies $\bbA_1$ is 
the following. This will be valuable in Section~\ref{S:ideals}. 

\begin{cor}\label{sametops}
The weak-$*$ and weak operator topologies coincide on $\flgee$. 
\end{cor}

\begin{rem}
Applied to the case of single vertex graphs, the proof of Theorem~\ref{a1thm} for $n=1$ 
provides a new proof of property $\bbA_1$ for the free semigroup algebras 
$\fL_n$, $n\geq 2$. In fact, this proof unifies the $H^\infty$ \cite{BFP} and $\fL_n$ 
cases. The previously known proof for $\fL_n$ \cite{DP2} of Davidson and Pitts 
relied on the existence of a pair of isometries with mutually orthogonal ranges in the commutant. 
The subclass of $\flgee$ which satisfy this extra condition are discussed further below. 
\end{rem}

\section{Partly Free Algebras}\label{S:partlyfree}

An interesting subclass of 
the algebras $\flgee$  discovered in \cite{KP2,KP1} are characterized by  commutant
$\flgee^\prime = \frgee$ containing a 
pair of isometries with mutually orthogonal ranges. In the terminology of 
\cite{KP2,KP1}, these are the $\flgee$ with {\it unitally partly free} 
commutant, as there is a unital injection of the free 
semigroup algebra $\fL_2$ into the commutant $\flgee^\prime$ in this case. 

We say that $G$ contains a {\it double cycle} if there are distinct cycles $w_i = kw_ik$,
$i=1,2$, of minimal length over some vertex $k$ in $G$. 
By a {\it proper infinite (directed) path} in $G$, we mean an infinite 
path $\omega =  e_{i_1}
e_{i_2}e_{i_3} \cdots$ in the edges of $G$ such that no edges are repeated 
and every finite segment corresponds to an allowable finite directed path in $G$. 
(Note that with our notation such an infinite path ends at the final vertex for $e_{i_1}$.) 
Say that $G$ has the {\it aperiodic path property} if  there 
exists an aperiodic infinite path; in other words,  there exists a 
proper infinite path or a double cycle in $G$.
Define the {\it attractor} at a vertex $k$ in
$G$ to be the set  consisting of $k$, together with all finite 
and infinite paths which end at $k$,
and all vertices that are initial vertices for paths ending at $k$.  
Then $G$ satisfies the
{\it uniform aperiodic path entrance property} if the attractor at every vertex 
includes an aperiodic infinite path.  

The following result is from \cite{KP1} 
(finite case) and \cite{KP2} (countably infinite case).

\begin{lem}\label{partlyfree}
The following assertions are equivalent for a countable directed 
graph $G$:
\begin{itemize}
\item[$(i)$] $G$ has the uniform aperiodic path entrance property.
\item[$(ii)$] $\flgee^\prime$ is unitally partly free.
\item[$(iii)$] $\flgee^\prime$ contains a pair of isometries $U$, $V$ with 
mutually orthogonal ranges; $U^* V =0$. 
\end{itemize}
\end{lem}

In \cite{KP2,KP1} things were  phrased differently. Recall that the commutant $\flgee^\prime =
\frgee$ is unitarily equivalent to $\fL_{G^t}$, 
where $G^t$ is the graph obtained from $G$  by reversing 
directions of all directed edges. Thus, as in \cite{KP2,KP1}, one could 
just as easily phrase this theorem in terms of $\flgee$ being partly free, proper infinite paths
defined by starting (instead of ending) at a vertex, and 
the graph $G$ satisfying a corresponding exit property. 

The $\flgee$ satisfying the condition in Lemma~\ref{partlyfree} 
form a large class of operator algebras \cite{KP2}.  
Indeed, the classification theorem from \cite{KP1}, which proved $G$ 
to be a complete unitary invariant of $\flgee$, shows that different 
graphs really do yield different algebras.  
Concerning the examples from \ref{examples}, observe that the graphs in 
$(ii)$ and $(v)$ satisfy the uniform aperiodic path entrance condition, 
hence $\flgee^\prime$ is unitally partly free in both cases. 
Further, notice that the algebra $\fL_{C_\infty}$ from $(iv)$ is unitally partly free since it
is unitarily equivalent to the commutant of $\fL_{C^t_\infty}$.


In connection with predual properties, the importance of  isometries with pairwise orthogonal 
ranges in the commutant of an algebra was realized in \cite{Berc} and \cite{DP3}. In particular,
the following result was proved in \cite{DP3}. 

\begin{lem}\label{prepredual}
If $\fA$ is a $\wot$-closed algebra which commutes with two isometries with orthogonal ranges,
then it has property $\bbA_{\aleph_0}$. 
\end{lem}

Thus 
we have the following consequence of Lemma~\ref{partlyfree}. 

\begin{cor}\label{predual}
If $G$ has the uniform aperiodic path entrance property, equivalently the 
commutant $\flgee^\prime$ is unitally partly free, then 
$\flgee$ has property $\bbA_{\aleph_0}$.
\end{cor}

There is a strong partial converse of this result, to which we now turn for the rest of this section. 

\begin{defn}
If $G$ is a directed graph for which the attractor at each vertex includes an infinite
directed path, then we say $G$ has the {\it uniform infinite path entrance property}. 
\end{defn}

The class of algebras $\flgee$ generated by the graphs which satisfy this property includes, for
example, all algebras satisfying the
conditions of Lemma~\ref{partlyfree}. But function algebras such as $H^\infty$ and the cycle
algebras $\fL_{C_n}$, $1 \leq n < \infty$, are also included. The following result shows that the
algebras of Lemma~\ref{partlyfree} truly stand apart in this class. They are the only algebras
which satisfy the $\bbA_{\aleph_0}$ property.

\begin{thm}\label{dichotomy}
Let $G$ satisfy the uniform infinite path entrance property. Then the following assertions are
equivalent: 
\begin{itemize}
\item[$(i)$] $G$ has the uniform aperiodic path entrance property.
\item[$(ii)$] $\flgee^\prime$ is unitally partly free.
\item[$(iii)$] $\flgee^\prime$ contains a pair of isometries $U$, $V$ with 
mutually orthogonal ranges; $U^* V =0$. 
\item[$(iv)$] $\flgee$ satisfies property $\bbA_{\aleph_0}$. 
\end{itemize}
\end{thm}

We first show that the cycle algebras do not satisfy $\bbA_{\aleph_0}$. 

\begin{lem}\label{cna2}
Let $1 \leq n < \infty$. The cycle algebra $\fL_{C_n}$ does not satisfy property $\bbA_{n+1}$. 
\end{lem}

\Prf
Suppose $\fL_{C_n}$ satisfies $\bbA_{n+1}$.
Choose a vertex $k$ in $C_n$ and put $\H = \H_{C_n}$. Then, since $P_k$ belongs to
$\fL_{C_n}$, the compression algebra $\fA = P_k \fL_{C_n}P_k |_{P_k\H}$ has property
$\bbA_{n+1}$ as well. To see this, suppose $\varphi$ is a functional on $\fA$. For $A\in
\fL_{C_n}$, let $A_k = P_k A P_k |_{P_k\H}$ and define a functional on $\fL_{C_n}$ by 
$\td{\varphi} (A) \equiv \varphi (A_k)$. If $\td{\varphi}$ can be realized as a vector functional,
$\td{\varphi}(A) = (A\xi,\eta)$, then 
\[
\varphi(A_k) = \varphi (P_k A_k P_k |_{P_k\H}) = \td{\varphi} (P_k A_k P_k) = 
(A_k P_k \xi, P_k \eta). 
\]

Let $w\in\fngee$ be the cycle of minimal
length in $C_n$ with $w=kwk$. 
A consideration of Fourier expansions for elements of $\fL_{C_n}$ reveals that 
$P_k\fL_{C_n} P_k$ is the subalgebra of $\fL_{C_n}$ given by 
\[
P_k\fL_{C_n} P_k = \wot \!\! - \!\! \Alg \big\{ L_w , P_k \big\}. 
\]
Moreover, $V = P_k L_wP_k |_{P_k \H}= L_w|_{P_k \H} $ is  unitarily equivalent to the
canonical
unilateral shift operator of multiplicity $n$. Indeed, $P_k$ commutes with $L_w$ and $V$ is a
pure isometry with 
\[
\rank (I-VV^*) = \rank \big( P_k(I-L_wL_w^*)P_k\big) = n.
\]
(One multiplicity is picked up for each of the $n$ infinite stalks in the Fock space tree.)
Thus, $\fA =\wot \! - \! \Alg \{V, I_{P_k\H}\}$ does not have property $\bbA_{n+1}$ 
since the shift of multiplicity $n$ does not satisfy $\bbA_{n+1}$  \cite{BFP}. 
This contradiction shows that
$\fL_{C_n}$ does not satisfy property $\bbA_{n+1}$.
\bx
 
\vspace{0.1in}

{\noindent}{\it Proof of Theorem~\ref{dichotomy}.} 
It follows from
Lemma~\ref{partlyfree} and Corollary~\ref{predual} that it suffices to prove $(iv)\Rightarrow
(i)$. Thus, 
suppose  $\flgee$ satisfies $\bbA_{\aleph_0}$. Given a subset $\S \subseteq V(G)$ of vertices in
$G$,
let $P_\S$ be the projection in $\flgee$ defined by $P_\S = \sum_{k\in\S} 
P_k$. Observe that every compression
algebra $P_\S \flgee P_\S|_{P_\S \H}$ satisfies $\bbA_{\aleph_0}$ since $P_\S$ 
belongs to $\flgee$; every functional $\varphi$
on $P_\S \flgee P_\S|_{P_\S \H}$ can be extended to $\flgee$ by defining 
$\td{\varphi}(A)\equiv \varphi(P_\S
A P_\S|_{P_\S \H})$ as in the previous proof. 

Now suppose $(i)$ does not hold. Then there is a vertex $k\in V(G)$ such that the attractor at
$k$ contains no aperiodic infinite path. In other words, there are no paths leaving double cycles
for $k$, and there are no infinite non-overlapping directed paths which end at $k$. Hence, the
uniform infinite path entrance property tells us there is a path from a cycle into $k$.  
Furthermore, by moving backwards along the paths which enter $k$, we can find a cycle $C$
from which there is a path into $k$ with the following properties: The only edges that enter the
vertices in $C$ are the edges which make up the cycle, and the cardinality $n$ of the vertex set
is equal to the number of edges in the cycle. In other words, there are no multiple edges
between vertices in $C$ (using the fact that $(i)$ fails), and there are no directed paths from
vertices outside $C$ to vertices
inside $C$ (using the uniform infinite path entrance property). Let $\S$ be the collection of
vertices in $C$, so that $|\S| = n$. It follows from the
choice of $C$ that the algebra $P_\S \flgee P_\S|_{P_\S \H}$ is unitarily equivalent to
$\fL_{C_n}$. This gives a contradiction to Lemma~\ref{cna2}, and hence $G$ must in fact
satisfy $(i)$. 
\bx

\begin{rem}
We wonder whether the initial restriction to $G$ satisfying the uniform infinite path entrance
property
is really necessary in Theorem~\ref{dichotomy}. It may simply be a convenient technical
assumption. It seems plausible to us that Theorem~\ref{dichotomy} could hold without making
this initial restriction on $G$. We would also expect that the equivalent conditions in this result
could be extended to include the related factorization properties  
$\bbA_n(n^2)$ and
$X_{0,1}$ \cite{Berc,DP3}, which involve norm control over the vectors 
chosen, and perhaps even property $\bbA_2$.  
\end{rem} 


 

\section{Ideals and Invariant Subspaces}\label{S:ideals}

In this section we give a detailed description of the $\wot$-closed ideal 
structure for $\flgee$.
The key ingredients in our analysis are the Beurling Theorem~\ref{beurlingthm} for $\flgee$
and the $\bbA_1$ property for functionals established in
Theorem~\ref{a1thm}. 
In \cite{DP2}, Davidson and Pitts described the ideal structure for  free semigroup algebras
$\fL_n$. 
For the sake of continuity in the literature, our presentation in this section will mirror
their approach whenever possible.  

\begin{defn}
Given a countable directed graph $G$, let $\Id_r(\flgee)$ and $\Id 
(\flgee)$ denote, respectively, the sets of all $\wot$-closed right and 
two-sided ideals. Note that $\Id_r(\flgee)$ and $\Id(\flgee)$ form 
complete lattices under the operations of intersections and $\wot$-closed 
sums. 
\end{defn}

To streamline the presentation, we shall let $\xi_\phi$ be the vector 
given by the following weighted sum of the vacuum vectors;  
$\xi_\phi = \sum_{k\in V(G)} \frac{1}{k}\, \xi_k$. This sum is finite 
precisely when $G$ has finitely many vertices. Observe that for all words 
$w\in \bbF^+(G)$ we have $L_w \xi_\phi = \frac{1}{k} \xi_w$, where $w = 
wk$. 

If $\J$ belongs to $\Id_r (\flgee)$, then the subspace 
\[
\ol{\J  \xi_\phi} = \ol{\J \flgee \xi_\phi} = \ol{\J \H_G} = \ol{\J \frgee \xi_\phi} = \frgee \ol{\J 
\xi_\phi}.
\]
Thus, the range subspace of $\J$ satisfies $ \ol{\J \H_G} = \ol{\J  \xi_\phi}$ and is
$\frgee$-invariant. If, in addition, $\J$ belongs to $\Id(\flgee)$, then $\flgee \J \xi_\phi = \J
\xi_\phi$, and $\ol{\J  \xi_\phi}$ is $\flgee$-invariant. Hence $\ol{\J  \xi_\phi}$ belongs to both
$\Lat (\flgee)$ and $\Lat(\frgee)$ when $\J$ is a two-sided ideal. 

On the other hand, if $\M$ belongs to $\Lat (\frgee)$, it follows that the set $\{A\in\flgee :
A\xi_\phi \in\M\}$ is contained in $\Id_r (\flgee)$. Indeed, this set is clearly $\wot$-closed and
for $X \in \flgee$, 
\[
AX \xi_\phi \in \ol{A\H_G} = \ol{A\frgee \xi_\phi} = \ol{\frgee A \xi_\phi} \subseteq \M. 
\]
Furthermore, if $\M$ is also $\flgee$-invariant, then this set evidently forms a two-sided ideal. 

The following theorem shows that ideals in $\Id_r (\flgee)$ and $\Id(\flgee)$ are fully
determined by their ranges.  

\begin{thm}\label{latticeiso}
Let $\mu : \Id_r(\flgee) \rightarrow \Lat (\frgee)$ be defined by $\mu (\J) = \ol{\J \xi_\phi}$.
Then $\mu$ is a complete lattice isomorphism. The restriction of $\mu$ to the set $\Id(\flgee)$
is a complete lattice isomorphism onto $\Lat(\flgee)\cap \Lat(\frgee)$. The inverse map 
$\iota$ sends a subspace $\M$ to 
\[
\iota (\M) = \big\{ A\in\flgee \,\, : \,\, A\xi_\phi \in \M \big\}. 
\]
\end{thm}

\Prf
We have observed above that the maps $\mu$ and $\iota$ map into the correct subspace lattice
and ideal lattice respectively. 

We show first that $\mu\iota$ is the identity map. Let $\M$ belong to $\Lat (\frgee)$. It is clear
that $\mu\iota (\M)$ is contained in $\M$. Conversely, let $\{ \xi_{k,j}\}_{k,j}$ be an
orthonormal basis for the $\frgee$-wandering subspace 
$\W = \M \ominus \sum_e \oplus R_e \M$, with $Q_k \xi_{k,j} = \xi_{k,j}$. Then from the
Beurling Theorem (the $\frgee$ version) we have  
\[
\M = \sum_{j,k} \oplus \frgee [\xi_{k,j}] = \sum_{j,k} \oplus \ran (L_{\xi_{k,j}}).
\]
As $L_{\xi_{k,j}}^* L_{\xi_{k,j}} = P_k$, it follows that $L_{\xi_{k,j}} \xi_\phi = \frac{1}{k}
\, \xi_{k,j}$
is in $\M$ and $L_{\xi_{k,j}}$ belongs to $\iota(\M)$. Hence 
\[
\M = \sum_{j,k} \oplus \ran (L_{\xi_{k,j}}) \subseteq \ol{\iota (\M) \H_G} = \ol{\iota (\M)
\xi_\phi} = \mu(\iota(\M)). 
\]
Thus, $\mu \iota (\M) = \M$, as required. 

To see that $\iota \mu$ is the identity, fix $\J$ in $\Id_r(\flgee)$ and let $\M = \mu(\J)$. It is
clear from the definitions that $\J$ is contained in $\iota \mu (\J)$. We first
show that for every $\xi$ in $\H_G$, 
\begin{eqnarray}\label{muiota}
\ol{\J \xi} = \ol{\iota \mu(\J) \xi}. 
\end{eqnarray}
From the Beurling Theorem, the cyclic $\flgee$-invariant subspace $\flgee [\xi]$ decomposes as
$\flgee[\xi] = \sum_{k\in \S} \oplus R_{\zeta_k}\H_G$, where $\S = \{ k\in V(G): P_k \xi \neq
0\}$, the vectors $\zeta_k = P_k \zeta_k$ are $\flgee$-wandering, and the $R_{\zeta_k}$ are
partial isometries in $\frgee$ with pairwise orthogonal ranges.  It follows that 
\[
\ol{\J \xi} = \ol{\J \flgee \xi} = \ol{\J \sum_{k\in\S} \oplus R_{\zeta_k}\H_G} = 
 \ol{ \sum_{k\in\S} \oplus R_{\zeta_k}\J \H_G} =  \sum_{k\in\S} \oplus R_{\zeta_k}\M. 
\]
But since $\mu \iota$ is the identity we have  $\mu(\iota\mu (\J)) = \mu(\J) = \M$, hence the
same computation for
$\ol{\iota\mu(\J) \xi}$ yields the same result. This yields (\ref{muiota}). 

Next suppose that $\varphi$ is a $\wot$-continuous linear functional on $\flgee$ which
annihilates the ideal $\J$. By Theorem~\ref{a1thm}, there are vectors $\xi, \eta$ in $\H_G$ with 
$\varphi (A) = (A\xi, \eta)$ for $A$ in $\flgee$. Since $\varphi(\J) = 0$, the vector $\eta$ is
orthogonal to $\ol{\J\xi} = \ol{\iota\mu(\J)\xi}$. Hence $\varphi$ also annihilates $\iota
\mu(\J)$. Thus, since the weak-$*$ and $\wot$ topologies on $\flgee$ coincide, we have $\J =
\iota \mu (\J)$ by the Hahn-Banach Theorem. 

Therefore we have shown that $\mu$ establishes a bijective correspondence between elements
of $\Id_r(\flgee)$ and $\Lat (\frgee)$ which maps $\Id(\flgee)$ onto $\Lat(\flgee)\cap
\Lat(\frgee)$, and that $\iota = \mu^{-1}$. It is elementary to verify that $\mu$ and $\iota$ are
complete lattice isomorphisms, and, since this may be accomplished in exactly the same way as
for the free semigroup algebras $\fL_n$ \cite{DP2}, we leave these remaining details to the
interested reader. 
\bx

Let us apply the theorem to the case of singly generated ideals. 

\begin{cor}\label{singly}
Let $A$ belong to $\flgee$. Then the $\wot$-closed two-sided ideal generated by $A$ is given
by $\big\{ X\in\flgee \, : \, X\xi_\phi \in \ol{\flgee A \H_G} \big\}$. 
\end{cor}

An interesting  special case occurs when $A = L_w$ is a partial isometry coming from a word
$w$ in $\bbF^+(G)$. The corresponding ideal is easily described in terms of Fourier expansions.

\begin{cor}
For $w$ in $\bbF^+(G)$, the $\wot$-closed two-sided ideal  generated by $L_w$ is given by 
\[
\Big\{ X \in\flgee \, : \, (X\xi_\phi, \xi_v) \neq 0 \qonlyif v=u_1wu_2; \, u_1,u_2\in\bbF^+(G)
\Big\}.
\]
\end{cor}



The theorem also leads to a simple characterization of the $\wot$-closure of the commutator
ideal of $\flgee$. Define the {\it $G$-symmetric Fock space} to be the subspace of $\H_G$
spanned by the vectors $\sum_{\sigma \in S_r} \xi_{\sigma(w)}$, where $w\in\fgeeplus$ with
$|w|=r$, $S_r$ is the symmetric group on $r$ letters, and $\sigma(w)$ is the word with letters in
$w$ permuted by $\sigma$. We put $\xi_{\sigma(w)} = 0$ and 
$L_{\sigma(w)} =0$ when $\sigma(w)$ is not an
allowable finite path in $G$. The terminology here is motivated by the case of a trivial graph
with one vertex and a number of loop edges. For a general graph though, there may be very little
`symmetry' associated with $\H_G$ (see Example~\ref{commutatoreg} below). 

By considering Fourier expansions it is not hard to see that the 
linear span of the commutators of 
the form $[L_v,L_w] = L_v L_w - L_w L_v$ for $v,w$ in $\fgeeplus$ is 
$\wot$-dense in the $\wot$-closure of the commutator ideal of $\flgee$. 
But a simple exercise in combinatorics shows that each of these commutators is determined by
elementary commutators. For instance, given edges $e,f,g$ in 
$E(G)$ observe that
\[
[L_{ef},L_g] = L_e [L_f,L_g] + [L_e,L_g]L_f. 
\]
It follows that the $\wot$-closure of the commutator ideal and its range 
subspace have the following form.

\begin{cor}\label{commutator} 
The $\wot$-closure of the commutator ideal
of $\flgee$ is the $\wot$-closed two-sided ideal given by 
\begin{eqnarray*}
\ol{\fC} &=& \Big< [L_v, L_w] : v,w \in\fgeeplus \Big> \\ 
&=& \Big< [L_e,L_f], [L_e,P_k] : e,f\in E(G), k\in V(G) \Big> 
\end{eqnarray*} 
The corresponding range subspace in $\Lat(\flgee) \cap
\Lat (\frgee)$ is 
\[ 
\mu(\ol{\fC}) = \big( \H_G^s \big)^\perp = \spn
\big\{ \xi_{uefv} - \xi_{ufev} : e\neq f , \, u,v\in\fgeeplus \big\}.  
\]
\end{cor}

Notice that an elementary commutator will typically collapse; for instance, $L_e L_f = 0$ if the
final vertex of $f$ is different than the initial vertex of $e$. In fact, $L_e L_f$ and $L_f L_e$ 
are both non-zero precisely when $ef$ (and $fe$) forms a cycle in $G$. 
The following class of examples differ greatly from the algebras $\fL_n$ 
and $H^\infty$. 

\begin{eg}\label{commutatoreg}
Let $G$ be a directed graph with no cycles; that is, no paths with the same initial and final
vertices. Then every commutator $[L_e, L_f] = L_e L_f - L_f L_e $ with $e,f\in E(G)$ is equal
to $L_eL_f$, $-L_fL_e$, or $0$. Further, since there are no loop edges, 
\[
[L_e, P_k] = \left\{ \begin{array}{cl}
L_e & \mbox{if $e=ek$} \\ 
-L_e & \mbox{if $e=ke$} \\ 
0  & \mbox{otherwise} 
\end{array}\right. 
\]
Thus, by Corollary~\ref{commutator}, the $\wot$-closure of the commutator 
ideal is equal to the
$\wot$-closed ideal generated by the $L_e$, 
\[
\ol{\fC} = \Big< L_e : e\in E(G) \Big>, 
\]
and the range is given by, 
\begin{eqnarray*}
\mu(\ol{\fC}) = \big( \H_G^s \big)^\perp &=& \spn \big\{ \xi_w : w\in\fgeeplus\setminus V(G)
\big\} \\
&=& \H_G \ominus \spn \{ \xi_k : k\in V(G) \}. 
\end{eqnarray*}
It also follows that $\flgee / \ol{\fC}$ is completely isometrically isomorphic to the algebra
$\spn\{ P_k : k \in V(G) \}$ in this case (see Section~\ref{S:distance}). 
\end{eg}

We finish this section with an investigation into factorization in right ideals. 

\begin{lem}\label{factorlemma}
Let $\{ L_{\zeta_j} : 1\leq j \leq s\}$ be a finite set of partial isometries in $\flgee$ with
pairwise orthogonal ranges $\M_j$. Let $\M = \sum_{j=1}^s \oplus \M_j$ and $\J = \iota (\M)$.
Then $\J = \{ A\in\flgee : \ran(A) \subseteq \M \}$ and every element $A$ of $\J$ factors as 
\[
A = \sum_{j=1}^s L_{\zeta_j} A_j \qwith A_j \in \flgee. 
\]
In particular, the algebraic right ideal generated by the finite set $\{ L_{\zeta_j} : 1\leq j \leq s\}$
coincides with $\J$. 

In the case of a countably infinite set of partial isometries $\{ L_{\zeta_j} : j
\geq 1\}$ in $\flgee$ with pairwise orthogonal ranges, every element $A$ of $\J$ (which is not a 
finitely generated algebraic ideal in this case) factors as a $\wot$-convergent sum 
\[
A = \wot-\sum_{j\geq 1} L_{\zeta_j} A_j \qwith A_j \in \flgee. 
\]
\end{lem} 

\Prf
We shall focus on the case where the $L_{\zeta_j}$ form a finite set. The countable case is
easily obtained from this analysis. 

Each $\M_j$ is $\frgee$-invariant, hence so is $\M$. Thus, if $A$ in $\flgee$ satisfies $A
\xi_\phi \in \M$, then $\ran(A) = \ol{A\H_G}$ is contained in $\M$, and we have 
\[
\J = \big\{ A\in \flgee \, : \, \ran(A) \subseteq \M \big\}. 
\]
This shows that $\J$ is a $\wot$-closed right ideal containing $\{ L_{\zeta_j} : 1\leq j \leq s\}$. 

On the other hand, since the projection onto $\M$ is given by $P_\M = \sum_{j=1}^s
L_{\zeta_j}L_{\zeta_j}^*$, for $A$ in $\J$ we have 
\[
A = \Big(  \sum_{j=1}^s L_{\zeta_j}L_{\zeta_j}^* \Big) A =  \sum_{j=1}^s L_{\zeta_j} A_j, 
\]
where $A_j =  L_{\zeta_j}^* A$. We finish the proof by showing that each $A_j$ is in $\flgee$.
First note that $A_j$ clearly commutes with the projections $Q_k\in\frgee = \flgee^\prime$. 
Further, as  $\frgee$-wandering vectors for $\M$, each $\zeta_j$ is orthogonal to 
$\sum_e \oplus R_e (\ran(A)) \subseteq \sum_e \oplus R_e\M$. Thus for $w\in \bbF^+(G)$, we
have $(R_e^*A^* \zeta_j, \xi_w ) = ( \zeta_j, R_e A \xi_w ) = 0$, and hence $R_e^*A^* \zeta_j 
=0$. 
Whence, using [Lemma~11.1, \cite{KP1}], given $e=ek$ in $E(G)$ we have  
\begin{eqnarray*}
A_jR_e - R_e A_j &=& L_{\zeta_j}^* A R_e - R_e L_{\zeta_j}^* A = (L_{\zeta_j}^* R_e -
R_e L_{\zeta_j}^* )A \\ 
&=& k^{2} \Big( \xi_\phi (R_e^* L_{\zeta_j} \xi_\phi )^*\Big) A = k \Big( \xi_\phi (A^* R_e^*
\zeta_j  )^*\Big) \\
&=& k \Big( \xi_\phi ( R_e^*A^* \zeta_j  )^*\Big) = 0.
\end{eqnarray*}
Therefore $A_j$ belongs to $\frgee^\prime = \flgee$, as required.
\bx

A special case of the lemma concerns the two-sided ideals in $\flgee$ generated by the partial
isometries from paths of a given length, $\{ L_w : w\in\fgeeplus, |w| =s\}$. 

\begin{cor}
For $s \geq 1$, every $A$ in $\flgee$ can be written as a sum 
\[
A = \sum_{|w|< s} a_w L_w + \sum_{|w| = s} L_w A_w, 
\]
where $a_w \in \bbC$ and $A_w \in \flgee$ for $w\in \fgeeplus$. In the case that there are
infinitely many paths of a given length, $\{ a_w \}_{|w| < s}$ belongs to $\ell^2$ and the sums
are $\wot$-convergent. 
\end{cor} 

We mention that the previous two results applied to the $\fL_n$ case include a uniqueness of
factorization. For the general $\flgee$ case this uniqueness does not hold; ostensibly because the
generators here are partial isometries instead of isometries. However, the elements $A_j$
in $\flgee$ from Lemma~\ref{factorlemma} can be chosen uniquely under the extra constraint
$L_{\zeta_j}^* L_{\zeta_j} A_j = A_j$. 

Together with
Theorem~\ref{latticeiso}, the previous lemma may be used to describe precisely when a right
ideal is finitely generated. The following result generalizes [Theorem~2.10, \cite{DP2}]. 

\begin{thm}\label{algebraicideal}
Let $\J$ be a $\wot$-closed right ideal in $\flgee$. Let $\M = \mu(\J)$ in $\Lat (\frgee)$ and let
$\W$ be the $\frgee$-wandering subspace for $\M$. If the sum of the dimensions of the
wandering subspaces $\{ Q_k \W : k\in V(G)\}$ is finite, $s< \infty$, then $\J$ is generated by
$s$ partial isometries with pairwise orthogonal ranges as an algebraic right ideal. When this
wandering dimension is infinite, $s = \infty$, $\J$ is not finitely generated as a $\wot$-closed
right ideal, but it is generated by countably many partial isometries as a $\wot$-closed right
ideal. 
\end{thm}

\begin{rem}
There is no analogous structure for the $\wot$-closed left ideals of $\flgee$. While some partial
results go through, there are factorization pathologies  in left ideals. Indeed, this was discovered
in \cite{Kfactor} for the case of  free semigroup algebras $\fL_n$. For instance, 
the algebraic left ideal determined by a generator of $\fL_n$ is not even norm closed. The basic
point is that these generators have proper factorizations inside the algebra, a property which is
exclusive to the noncommutative setting. 
\end{rem}

We conclude this section with a  comment on related work in the literature. 

\begin{rem}
The results on ideals of $\flgee$ presented in this section  generalize the characterization of
ideals in $H^\infty$ \cite{Douglas_text,Hoffman_text}, as well as work of Davidson and Pitts
\cite{DP1} in the case of free semigroup algebras $\fL_n$. In their work on quiver algebras,
Muhly and Solel \cite{MS1} briefly considered the ideal structure for $\flgee$, obtaining a
version of Theorem~\ref{latticeiso} in the special case that $G$ is a finite graph and satisfies a
certain {\it entrance condition}. We mention that their condition can be seen to be equivalent to
the finite graph case of the  {\it uniform aperiodic path entrance property} from
Theorem~\ref{partlyfree}. Thus, by Theorem~\ref{partlyfree},  Muhly and Solel actually 
established the lattice isomorphism theorem in the case that $G$ is finite and $\flgee$ has
unitally partly free commutant, and hence our result is an extension of theirs to the general
$\flgee$ case when $G$ has no sinks.  
\end{rem}

\section{Distance Formula to Ideals}\label{S:distance}

Let $\M_n(\flgee)$ denote the algebra of $n\times n$ matrices with entries in $\flgee$, equipped
with the operator norm in $\B(\H_G)^{(n)}$.  We prove the following distance formula to 
ideals in $\M_n(\flgee)$.

\begin{thm}\label{distformula}
Let $\J$ be a $\wot$-closed right ideal of $\flgee$, and let $\M= \ol{\J \H_G} = \ol{\J \xi_\phi}$
be the range subspace in $\Lat (\frgee)$.  Then for every $A$ in 
$\M_n(\flgee)$, we have 
$$
\dist (A, \M_n(\J)) = \|(P_\M^\bot \otimes I_n)A\|. 
$$
\end{thm}

\begin{rem}
Our proof is based on McCullough's distance formula for ideals in dual algebras \cite{McCull}. 
We shall proceed by establishing a variant of this result (Lemma~\ref{techlemma}), then
combine it with the Wold Decomposition Lemma~\ref{purewold} to obtain the distance
formula.  Note that Lemma~\ref{techlemma} does not use the fact that $\flgee$ is
$\wot$-closed, only that it is weak-$\ast$ closed.  One could prove the distance formula using
the
results of Section 3 (thus using the Wold decomposition implicitly, and avoiding
Lemma~\ref{techlemma}), but we have arranged the proof this way to emphasize the role of the
Wold decomposition.  The proof we give reduces to that given by McCullough for $H^\infty$,
and is similar in flavour to the Arias-Popescu proof  in the case of  free semigroup algebras
$\fL_n$ \cite{AP2}. 
\end{rem}

We first fix some notation. Let $1 \leq n \leq \aleph_0$.   Throughout this section, $Z$ will
denote a positive trace class
operator in $\B(\H_G^{(n)})$, factored as $Z=\sum_{i\geq 1} z_i  z_i^*$, with each $z_i$ in
$\H_G^{(n)}$.  Let $\tilde{z} = [z_1 \, z_2 \, \cdots ]^t$ denote the corresponding vector in the
infinite direct sum $\H_G^{(\infty)}$.  Note that $\|\tilde{z}\|^2 = \sum_{i\geq 1} \|z_i\|^2 = \tr
(Z)$.  We will let
$\fA=\M_n(\flgee)$, and for a $\wot$-closed right ideal $\J$ in $\flgee$, we let
$\fA_\J=\M_n(\J)$.

For a given positive trace class operator $Z$ on $\H_G^{(n)}$, let
\[
  \M(Z) = \ol{\fA^{(\infty)}\tilde{z}}  
\qand
  \N(Z) = \ol{\fA_\J^{(\infty)}\tilde{z}},
\]
and let $P_{\M(Z)}$ and
$P_{\N(Z)}$ be the projections onto $\M(Z)$ and $\N(Z)$ respectively.

\begin{lem}\label{techlemma}
  For $A\in \fA$,
$$
\dist(A,\fA_\J)=\sup_Z \|P_{\N(Z)}^\bot A^{(\infty)} P_{\M(Z)}\|
$$
where the supremum is taken over all positive trace class operators $Z$ in $\B (\H_G^{(n)})$.  
\end{lem}

\Prf
We must show that
$$
\inf_{B-A \in \fA_\J} \|B\| =\sup_Z \|P_{\N(Z)}^\bot A^{(\infty)} P_{\M(Z)}\|;
$$
the left-hand side being the definition of $\dist(A,\fA_\J)$.  We first prove the distance is
bounded above by  this supremum.

Since $\fA$ is weak-* closed, it is isometrically isomorphic to the dual 
of a Banach space; namely
$$
\fA_\ast = \C_1(\H_G^{(n)})\,\, / \,\, \fA_\bot
$$
where $\fA_\bot$ consists of the trace-class operators $T\in \C_1(\H_G^n)$ 
such that $tr(AT)=0$ for every $A\in \fA$.  Let $\pi$ denote the quotient 
map
$$
\pi: \C_1(\H_G^{(n)}) \longrightarrow \fA_\ast.
$$
Then the action of an element $A\in \fA$ on the predual $\fA_\ast$ is 
given by
$$
L_A(F)\equiv tr(A\tilde{F})
$$
where $\tilde{F}$ is any representative of the coset $F\in \fA_\ast$.  
Let $\fA_{\J\bot}$ denote the preannihilator of $\fA_\J$ in the predual $\fA_*$; that is, 
$$
\fA_{\J\bot} =\big\{ F\in\fA_* : L_A(F) = 0 \,\text{ for all }\, A\in 
\fA_\J \big\}.
$$
Given $F\in\fA_*$ and $\epsilon >0$, there exists a trace class operator $Y$ on $\H_G^{(n)}$
such that $\pi(Y) =F$, and $\|Y\|_1 < \|F\| + \epsilon$.  Let $Y^* =VZ$ be the polar
decomposition of $Y^*$, with $V$ a
partial isometry and $Z=(YY^*)^{1/2}$.
For this $Z$, let $\tilde{z}$ be as above.  Then given  $A\in \fA$, 
\begin{eqnarray*}
 ( A^{(\infty)}\tilde{z}, V^{*(\infty)}\tilde{z}) = \tr(AZV^*) 
  = \tr(AY) 
  =  L_A(F).
\end{eqnarray*}
Thus, when $A\in \fA_\J$, we have  
$
( A^{(\infty)}\tilde{z}, V^{*(\infty)}\tilde{z}) =0.
$

Now, since $\fA_\J^{(\infty)}\tilde{z} $ is dense in $\N(Z)$, it follows that
$$
\ran \big( P_{\M(Z)}V^{*(\infty)}P_{\M(Z)}\big) \subseteq \M(Z)\ominus \N(Z).
$$
We now compute, for any $A\in \fA_\J$,
\begin{eqnarray*}
  L_A(F) &=& ( A^{(\infty)}\tilde{z}, V^{*(\infty)}\tilde{z}) \\
  &=& ( A^{(\infty)}\tilde{z}, P_{\M(Z)} V^{*(\infty)}P_{\M(Z)}\tilde{z} ) \\
  &=& ( A^{(\infty)}\tilde{z}, P_{\N(Z)}^\bot P_{\M(Z)} V^{*(\infty)} P_{\M(Z)}\tilde{z}) \\
  &= & ( P_{\N(Z)}^\bot A^{(\infty)}\tilde{z}, P_{\M(Z)} V^{*(\infty)} P_{\M(Z)}\tilde{z} )
\end{eqnarray*}
Thus by the Cauchy-Schwarz inequality,
$$
|L_A(F)| \leq \|P_{\N(Z)}^\bot A^{(\infty)} P_{\M(Z)}\| 
\|V^{*(\infty)}\tilde{z} \| \|\tilde{z}\|.
$$
As $V$ is a partial isometry, $\|V^{*(\infty)}\tilde{z}\| \leq \|\tilde{z}\|$ and thus 
$$
\|\tilde{z}\|^2 = \tr(Z) =\|Y\| < \|F\| +\epsilon,
$$
implies (since $\epsilon$ was arbitrary) that 
\begin{align*}
  |L_A(F)| &\leq \|P_{\N(Z)}^\bot A^{(\infty)} P_{\M(Z)}\| \|\tilde{z}\|^2 \\
  &\leq \left( \sup_Z \|P_{\N(Z)}^\bot A^{(\infty)} P_{\M(Z)}\| \right) 
\|F\|.
\end{align*}
By the Hahn-Banach theorem, the functional $L_A$ on $\fA_{\J\bot}$ extends to a functional
$L$ on all of $\fA_*$ with 
$$
\|L\| \leq \sup_Z \|P_{\N(Z)}^\bot A^{(\infty)} P_{\M(Z)}\|
$$
and $L(F) =L_A(F)$ for $F\in \fA_{\J\bot}$.

Since $(\fA_*)^* =\fA$, there exists an operator $B$ in $ \fA$ such that $\|B\|=\|A\|$ and
\begin{align*}
L(F) &= L_B(F) \qfor F\in \fA_*, \\
L_A(F) &= L_B(F) \qfor F\in \fA_{\J\bot}.
\end{align*}
In particular, this means that $B-A$ belongs to $ (\fA_{\J\bot})^\bot$, and hence to $ \fA_\J$.

Thus, we have shown that there exists $B\in \fA$ such that $B-A\in \fA_\J$ and
$$
\|B\| \leq \sup_Z \|P_{\N(Z)}^\bot A^{(\infty)} P_{\M(Z)}\|,
$$
and so 
$$
\dist(A , \fA_\J) = \inf_{B-A \in \fA_\J} \|B\| \leq \sup_Z \|P_{\N(Z)}^\bot A^{(\infty)}
P_{\M(Z)}\| . 
$$

To prove the reverse inequality, observe that if $B-A$ belongs to $\fA_\J$, then 
$$
P_{\N(Z)}^\bot(B-A)^{(\infty)}P_{\M(Z)} =0,
$$
as $\N(Z)$ is the range of $\fA_\J^{(\infty)}$ restricted to $\M(Z)$.  Thus
$$
P_{\N(Z)}^\bot B^{(\infty)} P_{\M(Z)} = P_{\N(Z)}^\bot A^{(\infty)} P_{\M(Z)}.
$$
Hence, given $B$ with $B-A$ in $\fA_\J$, we have 
\begin{align*}
  \|B\| = \|B^{(\infty)}\| &\geq \|B^{(\infty)}P_{\M(Z)}\| \\
  &\geq \|P_{\N(Z)}^\bot B^{(\infty)} P_{\M(Z)}\| \\
  &= \|P_{\N(Z)}^\bot A^{(\infty)} P_{\M(Z)} \|
\end{align*}
Therefore, 
$$
\inf_{B-A \in \fA_\J} \|B\| \geq \sup_Z \|P_{\N(Z)}^\bot A^{(\infty)} P_{\M(Z)}\|.
$$
This concludes the proof of the lemma. 
\bx

\vspace{0.1in}

{\noindent}{\it Proof of Theorem~\ref{distformula}.} 
Suppose $A, B\in \fA$ and $A-B\in \fA_\J$.  Then, since $P_\M^\bot \otimes I_n$ annihilates
$\fA_\J$, 
$$
(P_\M^\bot \otimes I_n) B = (P_\M^\bot \otimes I_n) A.
$$
Hence 
$
   \|B\| \geq \|(P_\M^\bot \otimes I_n) B \| 
   = \| (P_\M^\bot \otimes I_n) A \|, 
$
and we have  
$$
\dist(A, \fA_\J) = \inf_{B-A \in \fA_\J} \|B\| \geq \| (P_\M^\bot \otimes 
I_n) A \|
$$

To prove the reverse inequality, we employ the previous lemma and the Wold 
Decomposition.  By Lemma~\ref{techlemma}, it 
will suffice to show that 
$$
\|P_{\N(Z)}^\bot A^{(\infty)} P_{\M(Z)}\| \leq \| (P_\M^\bot \otimes I_n) A \|
$$
for each positive trace class operator $Z$.

For $T_1, \ldots, T_n$ in $\flgee$, let $\tilde{T}$ denote the $n$-tuple $\td{T} = (T_1,\dots , 
T_n)$, treated as a column
vector.  We begin by constructing an auxilliary Hilbert space $\K$ in the following manner: 
Given a positive trace class
operator $Z$, define a Hermitian form on the set of $n$-tuples of elements of $\flgee$ via
$$
( \tilde{T},\tilde{S}) = \sum_{i\geq 1} (\td{T}^t z_i, \td{S}^t 
z_i)_{(\H_G)^{(n)}} = \sum_{i,k,l}( T_k z_i^k, S_l z_i^l )_{\H_G},
$$
where $z_i^k$ is the $k$th coordinate from $z_i\in \H_G^{(n)}$. 
The collection of $n$-tuples $\tilde{T}$ equipped with this Hermitian form gives a pre-Hilbert
space; taking the quotient modulo null vectors and closing, we obtain a Hilbert space which we
denote $\H_G(Z)$.  We let $\K = (\H_G(Z))^{(n)}$ denote the direct sum of 
$n$ copies of $\H_G(Z)$. 

If $A=(A_{jk})_{j,k=1}^n$ belongs to  $\M_n(\flgee)$, we may define a 
bounded operator $M_{A,Z}$
on $\K$ by $M_{A,Z} \equiv \big( A_{jk}^{(n)} \big)_{j,k=1}^n$. The effect 
of $M_{A,Z}$ on a vector $(\td{T}_1, \ldots, \td{T}_n)$ in $\K$ will be 
$$
M_{A,Z} (\td{T}_1, \ldots, \td{T}_n) =\left( \sum_{k=1}^n 
A_{1k}(\td{T}_k)_k, \ldots, \sum_{k=1}^n A_{nk} (\td{T}_k)_k  
\right)
$$
Now, let $\td{B}_j$ denote the transpose of the $j$th
row of $B\in M_n(\flgee)$, and define a map $\Phi$ by
$$
\Phi(B^{(\infty)}\tilde{z})= (\td{B}_1, \ldots, \td{B}_n) 
$$
between dense subspaces of $\M(Z)$ and $\K$.  An elementary calculation shows that $\Phi$ is
a surjective isometry between these dense subspaces, and so extends to a unitary map from
$\M(Z)$ onto
$\K$.  It is clear that, if we let $\K_\J$ denote the closure in $\K$ of 
vectors for which each entry
of each summand is an operator in $\J$, then $\Phi$ carries $\N(Z)$ (the range of the resriction
of $\fA_\J^{(\infty)}$ to $\M(Z)$) onto $\K_\J$.  Moreover, for each $A\in 
\fA$, this map
intertwines $A^{(\infty)}|_{\M(Z)}$ and $M_{A,Z}$.  Hence these maps are unitarily
equivalent, 
$$
A^{(\infty)}|_{\M(Z)} \simeq M_{A,Z}.
$$

Consider now the operators $\Psi_e$ on $\H_G(Z)$, defined by 
\[
\Psi_e(T_1,\dots , T_n) \equiv (L_eT_1,\dots ,L_eT_n).  
\]
The set $\{\Psi_e : e\in E(G)\}$ satisfies the hypotheses of
Lemma~\ref{purewold} with each $\dim \Psi_e^* \Psi_e \W \leq n$, and so 
there exists a unitary operator 
$$
U:\H_G(Z) \longrightarrow \H \equiv \sum_{k\in \S} \oplus (Q_k 
\H_G)^{(\alpha_k)}
$$
such that 
$$
U\Psi_e U^* = \sum_{k\in\S}\oplus L_e^{(\alpha_k)}|_{(Q_k \H_G)^{(\alpha_k)}} \qfor e\in
E(G), 
$$
where $\S $ is the subset of vertices $V(G)$ determined as in Lemma~\ref{purewold} by the
non-zero subspaces amongst $\{ \Psi_e^* \Psi_e \W : e\in E(G) \}$.

The ampliation $U^{(n)}$ is a unitary map from $\K$ onto $\H^{(n)}$. 
Under $U$, the $(j,k)$ block
$A_{jk}^{(n)}$ of $M_{A,Z}$ is taken to $A_{jk}^{(n)}|_{\H}$.  We now let $M_A$ denote
the operator on $\H_G^{(n^2)}$ which has the form of an $n\times n$ block matrix whose
$(j,k)$ block is $A_{jk}^{(n)}$.  With these definitions, $U^{(n)}M_{A,Z} (U^{(n)})^*
=M_A|_{\H^{(n)}}$.  Moreover, recalling that $\M=\ran(\J) = \ol{\J \H_G} \subseteq \H_G$, we
see that $U^{(n)}$ maps $\K_\J$ onto
$\M^{(n^2)}\cap \H^{(n)}=(\M^{(n)}\cap \H)^{(n)}$.  Thus, it follows that 
\begin{eqnarray*}
  P_{\N(Z)}^\bot A^{(\infty)} P_{\M(Z)} \simeq P_{\K_\J}^\bot M_{A,Z} \simeq (P_{\M}^\bot
\otimes I_{n^2}) M_A |_{\H^{(n)}}
\end{eqnarray*}

Finally, let $P$ denote the projection of $\H_G^{(n)}$ onto $\H$.  By  
reordering the summands
of $\H_G^{(n^2)}$ with a canonical shuffle, $M_A$ is seen to be unitarily 
equivalent to the direct sum of $n$ copies of
$A$ with itself.  This reordering fixes $P_\M^\bot \otimes I_{n^2}$, and takes $P$ onto another
projection, say $P^\prime$.  Thus, for every $A\in \M_n(\flgee)$, 
\begin{align*}
  \|P_{\N(Z)}^\bot A^{(\infty)} P_{\M(Z)}\| &= \| (P_{\M}^\bot \otimes I_{n^2}) M_A
|_{\H^{(n)}} \| \\
  &= \| (P_{\M}^\bot \otimes I_{n})^{(n)} A^{(n)} P^\prime \| \\
  &\leq \| (P_\M^\bot \otimes I_n) A \|,
\end{align*}
and this completes the proof.
\bx

The conclusion of Theorem~\ref{distformula} for two-sided ideals follows from the
result for right-sided ideals in precisely the same way as the $\fL_n$ case
\cite{DP3}. In particular, the following result may be applied to the $\wot$-closure of the
commutator ideal of $\flgee$, as well as the radical of $\flgee$ \cite{KP1}. 

\begin{cor}\label{twosideddistanceformula}
Let $\J$ be a $\wot$-closed two-sided ideal in $\flgee$, with range subspace $\M = \mu(\J)$.
Then $\flgee / \J$ is completely isometrically isomorphic to the compression algebra
$P_\M^\perp \flgee P_\M^\perp = P_\M^\perp \flgee $. 
\end{cor}

\section{Carath\'{e}odory Theorem}\label{S:caratheodory}
     
As an application of the distance formula, in this section we prove the analogue of the
Carath\'{e}odory
Theorem \cite{AP2,DP3} for the algebras $\flgee$. 
The   Carath\'{e}odory problem  specifies an initial segment  of the general Fourier series for
elements of $\flgee$, then asks when the 
segment can be completed to an operator in $\flgee$ of norm at most
one. 

We say that a subset $\Lambda$ of $\fgeeplus$ is a {\it left lower set} if $u$ belongs to
$\Lambda$ whenever $w=uv$ does. In other words, $\Lambda$ is closed under taking left
subpaths; these are the subpaths of $w$ obtained by moving from a vertex in $w$ along the rest
of $w = yw$ to its final vertex $y$. 
Let $P_\Lambda$ be the orthogonal projection onto the subspace $\M_\Lambda 
= \spn \{ \xi_w : w\in \Lambda \}$. Observe that $\M_\Lambda$ is invariant 
for $\frgee^*$, and hence $\M_\Lambda^\perp$ belongs to $\Lat (\frgee)$. 
Recall that elements $A$ of $\flgee$ have Fourier expansions $A \sim 
\sum_{w\in \fgeeplus} a_w L_w$, and hence it 
follows that elements of the matrix algebras $\M_k (\flgee)$ have natural 
Fourier expansions as well. 

\begin{thm}\label{carath}
Let $G$ be a countable directed graph and suppose $\Lambda$ is a left 
lower set of $\fgeeplus$. Given $k\geq 1$, let $\{ C_w : w\in \Lambda \}$ 
be matrices in $\M_k(\bbC)$. Then there is an element $A$ in the unit ball 
of $\M_k (\flgee)$ with Fourier coefficients $A_w = C_w$ for $w\in 
\Lambda$ if and only if 
\[
\Big|\Big| (P_\Lambda \otimes I_k ) \Big( \sum_{w\in\Lambda} L_w \otimes 
C_w \Big) \Big|\Big| \leq 1.
\]
When $\Lambda$ is an infinite set this sum is understood to converge 
$\wot$ in the Cesaro sense. 
\end{thm}

\Prf
This result can be proved as an immediate consequence of the distance 
formula. The set of elements of $\flgee$ which interpolate the zero data 
$C_w \equiv 0$ for $w\in\Lambda$ is equal to the $\wot$-closed right ideal 
$\J$ with range $\mu (\J) = P_\Lambda^\perp \H_G$. Thus, if we are given 
data $C_w$, the desired element exists precisely when the distance from 
$\sum_{w\in\Lambda} L_w \otimes C_w$ to $\M_k (\J)$ is at most one, and 
the result follows from Theorem~\ref{distformula}. 
\bx

The simplest example of a lower set is the set $\fgeeplus_k$ of all paths 
in $G$ of length at most $k$. The (two-sided) ideal in this case is given 
by 
\[
\iota (\M^\perp_{\fgeeplus_k}) \equiv \flgee^{0,k} = \big\{ A\in \flgee : A = 
\sum_{|w|> k} a_w L_w \big\}. 
\]
Thus we obtain the following. Let $P_{\Lambda_k}$ be the projection onto the 
subspace $\M_{\fgeeplus_k}$. 

\begin{cor}
Given a formal power series $\sum_{|w|\leq k} a_w w$ in the 
semigroupoid algebra $\bbC \fgeeplus$, there is an element $A$ in $\flgee$ 
with $||A|| \leq 1$ and $A - \sum_{w\in \fgeeplus_k} a_w L_w$ in $\flgee$ 
if and only if 
\[
\Big|\Big| P_{\Lambda_k} \Big( \sum_{w\in\fgeeplus_k} a_wL_w\Big) \Big|\Big| \leq 1.
\]
\end{cor}

\begin{rem}
There is  a different, more self-contained proof of the Carath\'{e}odory Theorem for 
$\flgee$ which is worth discussing since it yields  general information on elements of 
$\flgee$. Furthermore, this alternative approach gives a new proof for 
free semigroup algebras $\fL_n$ which generalizes the Parrot's Lemma cum 
Toeplitz matrix approach for $H^\infty$. 
\end{rem}

We first establish some notation and make some simple observations. 
Given a lower set $\Lambda$ in $\fgeeplus$, let $\Lambda_k = \{ 
w\in\Lambda : |w|\leq k\}$ for $k\geq 0$ and let $E_k $ be the  projection
of $\H_G$ onto $\spn \{ \xi_w : w\in\Lambda_k \}.$ 
Then each $E_k$ is a subprojection of $P_\Lambda$ and $P_\Lambda = \sum_{k\geq 0} \oplus
F_k$, where $F_{k+1} = E_{k+1} - E_k$ and $F_0 = E_0$. 
The lower set property yields the identity  
\begin{eqnarray}
E_{k+1} R_e E_{k} = E_{k+1} R_e \qfor  e\in E(G). 
\end{eqnarray}
Also note
that $Q_y = \sum_{w=wy} \oplus \xi_w \xi_w^*$ for every vertex $y\in V(G)$, where $\xi_w
\xi_w^*$ is the rank one projection onto the span of $\xi_w$, and that 
\[
R_e( \xi_w \xi_w^*)R_e^* = (R_e \xi_w) (R_e \xi_w)^* = \xi_{we} \xi_{we}^* .
\] 
Further, it is clear that the projections $\{ Q_y, P_\Lambda, E_k\}$ are mutually commuting; for
instance, $Q_y P_\Lambda E_k$ is the projection onto the subspace 
\[
Q_y P_\Lambda E_k \H_G = \spn \{ \xi_w : w\in\Lambda_k , w=wy\}. 
\]
We need the following estimates for elements of $\flgee$ to apply Parrot's
Lemma. 

\begin{prop}\label{parrotsetup}
Let $G$ be a countable directed graph and let $\Lambda$ be a lower set in $\fgeeplus$. 
Given $X\in\flgee$, define 
\[
A_k = E_k X E_k \qand B_k = E_{k+1} X (E_{k+1} - E_0) \qfor k\geq 0. 
\]
Then 
\[
||B_k|| \leq ||A_k|| \leq ||X|| \qfor k\geq 0.
\] 
\end{prop}

\Prf
Fix $k\geq 0$. 
As $Q_y$ commutes with $X$ and each $E_k$, we have $A_k = \sum_y \oplus A_k Q_y$.
Thus, 
\[
||A_k || = \sup_y || E_k X Q_y E_k ||. 
\]
On the other hand, 
\[
B_k = E_{k+1} X (E_{k+1} - E_0) E_{k+1} 
= \sum_{w\in\Lambda_{k+1}\setminus\Lambda_0} E_{k+1} X (\xi_w \xi_w^*) E_{k+1}. 
\]
But this sum may be written as
\[
B_k = \sum_{e\in E(G)} \sum_{u\in\Lambda_k} E_{k+1} X  (\xi_{ue} \xi_{ue}^*) 
E_{k+1}
\]
Indeed, since $\Lambda$ is a lower set, every $w\in \Lambda_{k+1}\setminus \Lambda_0$ is of
the form $w=ue$ for some $u\in\Lambda_k$ and edge $e$. On the other hand, if $w=ue$, with
$u\in\Lambda_k$, is not in $\Lambda_{k+1}$, then $(\xi_w \xi_w^*) E_{k+1} = \xi_w
(E_{k+1} \xi_w)^* = 0$, so the corresponding term in this sum vanishes. Thus, we have 
\begin{eqnarray*}
B_k &=& \sum_{e\in E(G)} \sum_{w\in\Lambda_k} E_{k+1} X R_e (\xi_w \xi_w^*) R_e^*
E_{k+1} \\
&=& \sum_{e\in E(G)} \sum_{w\in\Lambda_k} E_{k+1} R_e \big( E_{k} X (\xi_w \xi_w^*)
E_{k}\big) R_e^*E_{k+1} \\
&=& \sum_{y\in V(G)} \sum_{e=ye}\oplus E_{k+1} R_e \left( \sum_{w\in\Lambda_k; \,
w=wy} E_{k} X
(\xi_w \xi_w^*) E_{k}\right) R_e^* E_{k+1} \\
&=& \sum_y \sum_{e=ye}\oplus  E_{k+1} R_e \big(  E_{k} X Q_y  E_{k}\big) R_e^*E_{k+1}
.
\end{eqnarray*}
Therefore, as the ranges of $\{ E_{k+1}R_e:e\in E(G)\}$ are pairwise orthogonal for fixed
$k\geq 0$,  it follows that 
\begin{eqnarray*}
||B_k|| &=& \sup_y \sup_{e=ye} ||E_{k+1} R_e (E_k X Q_y E_k) R_e^*E_{k+1} || \\ 
&\leq& \sup_y || E_k X Q_y E_k || =  ||A_k||,
\end{eqnarray*}
for $k\geq 0$, and this completes the proof. 
\bx

These are precisely the estimates required to apply Parrot's Lemma
in the `bottom left corner' argument for the proof of the 
Carath\'{e}odory Theorem. In the $H^\infty$ case, equality is achieved with 
$||B_k || = ||A_k||$. This is easily seen from the Toeplitz matrix perspective 
for elements of $H^\infty$. This is also the case for free semigroup 
algebras $\fL_n$. Let us discuss the special case of $\fL_2$ to illustrate 
this generalized Toeplitz matrix structure. 

\begin{eg}
Let $\{1,2\}$ be the (noncommuting) generators of $\bbF^+_2$. The 
full Fock space has an orthonormal basis given by 
$\{ \xi_\phi, \xi_w : w\in \bbF_2^+ \}$. Elements $X$ of $\fL_2$ have 
Fourier expansions given by their action on the vacuum vector; 
$X \xi_\phi = \sum_{w\in\bbF_2^+\cup\{\phi \}} a_w \xi_w$. Let $\Lambda = 
\bbF_2^+ \cup \{ \phi \}$ be the trivial lower set, so that $P_\Lambda = I$. Consider the ordering
for $\Lambda_2 = \{ w : |w|\leq 2\}$ given by 
\[
\Lambda_2 = \Big\{ \phi, \{1, 2\}, \{1^2, 12, 21, 2^2\}\Big\}. 
\]
The projection $E_2$  has range space $\spn \{ \xi_w : w\in\Lambda_2  \}$. 
Hence the compression of $X \sim \sum_w a_w L_w$ in $\fL_2$ to $E_2$ is unitarily
equivalent to 
\[
A_2 = E_2 X E_2 \simeq \left[
\begin{matrix}
\begin{array}{c|cc|cccc}
a_\phi & 0 & 0 & 0 & 0 &0& 0 \\
\hline
a_1 & a_\phi & 0 & 0 & 0 &0& 0 \\
a_2 & 0 & a_\phi  & 0 & 0 &0& 0 \\
\hline
a_{1^2} & a_1 & 0  & a_\phi & 0 &0& 0 \\
a_{12} & 0 & a_1  & 0 & a_\phi &0& 0 \\
a_{21} & a_2 & 0  & 0 & 0 & a_\phi & 0 \\
a_{2^2} & 0 & a_2  & 0 & 0 & 0 & a_\phi 
\end{array} 
\end{matrix}\right]
\]
On the other hand, $B_2 = E_3 X(E_3 - E_0)$ can be seen to be unitarily equivalent to
$A_2\otimes I_2$ where $I_2$ is the scalar $2\times 2$ identity matrix. Thus, in this case
$||B_2|| = ||A_2||$. More generally, if $\Lambda = \fnplus \cup \{\phi \}$ and $X$ belongs to
$\fL_n$, the operators $B_k$ are unitarily equivalent to $A_k \otimes I_n$. The `$n$-branching'
in $\fnplus$ creates $n$-ampliations in the generalized Toeplitz matrices. 
\end{eg}

\begin{rem}
We mention that when $G$ has loop edges, it is also possible to derive a version of Pick's 
interpolation theorem \cite{AP2,DP3,Sar1} for $\flgee$.  Indeed, in 
\cite{KP1} it was shown how loop edges over vertices explicitly define 
eigenvectors for $\flgee^*$, with corresponding eigenvalues making up  
unit balls in the complex spaces $\bbC^k$. But the notation is rather  
cumbersome, and thus, as it is a direct generalization of the Pick 
theorem from 
\cite{AP2,DP3} for $\fL_n$, we shall not present it here. 
\end{rem}

{\noindent}{\it Acknowledgements.}
We are grateful to members of the Department of Mathematics at Purdue University for kind
hospitality during preparation of this article. The first author was  partially 
supported by a Purdue VIGRE Post-doctoral  Fellowship.  
The second author was partially supported by a Canadian NSERC
Post-doctoral Fellowship.


%

\end{document}